\documentclass{article}

\newif\ifarxiv
\arxivtrue

\usepackage{graphicx}
\usepackage{amsmath,amssymb,amsthm,bm}
\usepackage[a4paper]{geometry}
\usepackage{algorithm,algorithmic}
\usepackage[colorlinks=false,pdfborder={0 0 0}]{hyperref}
\usepackage{authblk}
\usepackage{color}

\graphicspath{{fig/}}
\DeclareGraphicsExtensions{.mps,.pdf,.eps,.png}

\newenvironment{keywords}{\medskip\textbf{Keywords:}}{}
\newenvironment{AMS}{\medskip\textbf{AMS subject classifications (2020).}}{}

\newtheorem{theorem}{Theorem}

\newtheorem{remark}{Remark}
\newtheorem{lemma}{Lemma}

\newtheorem{corollary}{Corollary}
\theoremstyle{plain}

\DeclareMathOperator{\diag}{diag}

\newcommand*{\trans}{^{\top}}

\newcommand*{\abs}[1]{\lvert#1\rvert}

\newcommand*{\norm}[1]{\Vert#1\rVert}

\def\adots{\mathinner{\mkern2mu\raise1pt\hbox{.}\mkern2mu
    \raise4pt\hbox{.}\mkern2mu\raise7pt\hbox{.}\mkern1mu}}

\newcommand*{\macheps}{\bm u_{w}}
\newcommand*{\machepslow}{\bm{u}_{\mathrm{\ell}}}
\newcommand*{\machepsmedm}{\bm{u}_{\mathrm{m}}}

\newcommand*{\bigO}{O}

\DeclareMathOperator{\fllow}{f{}l_{16}}
\DeclareMathOperator{\flworking}{f{}l_{64}}

\newcommand{\Acal}{\mathcal{A}}
\newcommand{\Acalhat}{\widehat{\mathcal{A}}}

\newcommand{\CC}{\mathcal{C}}

\newcommand{\GG}{\mathcal{G}}
\newcommand{\GGhat}{\widehat{\mathcal{G}}}
\newcommand{\LL}{\mathcal{L}}
\newcommand{\LLhat}{\widehat{\mathcal{L}}}

\newcommand{\UU}{\mathcal{U}}
\newcommand{\UUhat}{\widehat{\mathcal{U}}}

\title{Adaptive accuracy-driven splitting in Ozaki matrix multiplication and application to block LU factorization}
\author[1]{Erin Carson}
\author[1]{Yuxin Ma}
\affil[1]{Department of Numerical Mathematics, Faculty of Mathematics and Physics, Charles University, Sokolovsk\'{a} 49/83, 186 75 Praha 8, Czechia}

\begin{document}
\maketitle

\begin{abstract}
The Ozaki scheme for matrix multiplication computes an accurate matrix product by splitting the input matrices and evaluating a sequence of low-precision matrix products.
Its accuracy and computational cost depend strongly on the number of splits.
We develop accuracy-driven criteria for selecting this number in both Ozaki matrix multiplication and in block LU factorization using Ozaki trailing updates.  
For matrix multiplication, we derive normwise and componentwise criteria that determine how many split products are required to achieve a given accuracy requirement.
For block LU, we present a fixed normwise strategy with essentially no selection cost, an LU-aware adaptive normwise strategy that can reduce the number of splits in subsequent trailing updates, and a componentwise strategy that uses information accumulated from preceding updates to guide the later selections.  
In the componentwise strategy, when the required number of splits exceeds a prescribed maximum, the corresponding update is instead carried out using ordinary matrix multiplication.  
We also provide rigorous backward error bounds for these approaches.  
Numerical experiments show that the adaptive
normwise strategy can reduce the number of splits while preserving normwise accuracy, whereas the componentwise strategy can provide smaller componentwise backward errors when stronger local accuracy is needed.
\medskip


\end{abstract}


\section{Introduction}
\label{sec:introduction}

Low-precision matrix multiplication has become substantially faster than conventional double-precision arithmetic on modern accelerator architectures.
A representative example is NVIDIA's Tensor Core hardware, which is designed to provide very high throughput for low-precision matrix products.
For example, the NVIDIA B200 has a quoted peak throughput of 4.5~Pflop/s for dense FP8 Tensor Core operations, compared with 37~Tflop/s for dense FP64 Tensor Core operations, a difference of more than a factor of 121~\cite{NVIDIABlackwell}.
Such a large throughput gap creates a strong incentive to reformulate double-precision numerical kernels so that most of their cost can be shifted to low-precision matrix multiplication while retaining double-precision accuracy.

The Ozaki scheme provides one such approach.
It was introduced by Ozaki, Ogita, Oishi, and Rump in 2012 as an accurate matrix multiplication method based on an error-free transformation of a matrix product~\cite{Ozaki2012}.
The input matrices are decomposed into a number of split matrices, and the target product is represented by a sum of several ordinary matrix products.
In the original formulation, these products are evaluated in the same arithmetic as the input data, with the splitting and accumulation designed so that the final result is significantly more accurate than a direct matrix multiplication.

A major development came in 2020, when Mukunoki, Ozaki, Ogita, and Imamura proposed a mixed precision version for Tensor Cores~\cite{Ozaki2020}.
In their method, the split matrices are stored in FP16, the split products are evaluated by Tensor Cores with FP32 accumulation, and the resulting products are rescaled and accumulated in FP64.
The performance gain therefore comes from replacing most of the FP64 matrix multiplication work by much faster low-precision Tensor Core products.
The same basic idea has subsequently been extended to integer matrix multiplication units, where the split products are evaluated by high-throughput integer arithmetic~\cite{OOY2024,OUI2026,UOI2025}.
More recent work has also revisited floating point realizations on newer low-precision hardware such as FP8 Tensor Cores~\cite{Mukunoki2026}.

In this work, however, we focus on the floating point mixed precision Ozaki matrix multiplication proposed in~\cite{Ozaki2020}.
Compared with an integer-based realization, this formulation has several practical advantages.
First, it can be built directly on standard mixed precision GEMM interfaces such as \texttt{cublasGemmEx}, without introducing an additional fixed-point representation for the split matrices.
Second, the floating point formulation naturally accommodates the scaling used in the splitting procedure and avoids the explicit integer encoding, range management, and integer-to-floating point recovery required by an integer implementation.
Third, the approach is directly applicable on a wide range of accelerators that provide fast low-precision floating point matrix multiplication.
This point is becoming increasingly relevant as recent accelerator generations place substantial emphasis on low-precision floating point formats such as FP8 in addition to, or instead of, further increasing INT8 capability~\cite{Mukunoki2026}.
Thus, although integer Ozaki schemes can offer superior performance on suitable hardware, the floating point formulation remains an attractive and portable basis for studying how the Ozaki scheme can be incorporated into higher-level numerical algorithms.

For the Ozaki matrix multiplication, the number of splits, denoted by $d$, directly affects both the computational cost and the accuracy of the result.
A larger $d$ generally improves the accuracy, but requires more split matrix products and therefore reduces performance.
Conversely, a smaller $d$ can improve performance, but may not provide sufficient accuracy.
It is therefore important to choose $d$ according to the accuracy actually required by the computation.
In this work, we study how to select an appropriate $d$ for both the Ozaki matrix multiplication itself and block LU factorization using Ozaki matrix multiplication for the trailing matrix updates.
Related recent work~\cite{SABB2026} develops automatic dynamic precision based on an exponent span capacity estimator to determine the number of splits needed for Ozaki-based FP64 emulation on reduced precision Tensor Cores, together with automatic fallback to native FP64 when appropriate.
Very closely related to our matrix-multiplication analysis, \cite{ADFMT2025} derives componentwise error bounds for an integer-arithmetic realization of the Ozaki scheme and proposes a data-dependent strategy that allows different numbers of slices for the two input matrices to reduce computational cost, whereas~\cite{UOI2026-error} develops a deterministic error analysis for the distinct CRT-based Ozaki-II scheme and uses it to estimate the number of low-precision matrix multiplications required for a desired accuracy.

Our focus is the floating point Ozaki scheme of~\cite{Ozaki2020}.
We derive explicit normwise and componentwise criteria for selecting \(d\) according to a prescribed accuracy requirement, and further develop application-aware choices of \(d\) for the trailing updates in block LU factorization through a backward error analysis.
For normwise stability, we consider both a fixed choice of \(d\), which introduces essentially no additional selection cost, and an adaptive strategy based on prefix row and column norms of the computed factors.
For componentwise stability, the first trailing update uses the componentwise matrix-multiplication selector to initialize a nonnegative reference matrix \(\GG\), which is then used to select \(d\) in the subsequent updates; if the required number of splits exceeds a prescribed maximum, the corresponding update is evaluated by ordinary matrix multiplication.
We further derive backward error bounds for all of these strategies.
Unlike~\cite{ADFMT2025}, our componentwise selection criterion uses a computable lower bound that accounts for the contributions of corresponding rows and columns to the matrix product, rather than relying on separate componentwise relative truncation bounds for the two input matrices.
We further develop adaptive selection strategies tailored to the global backward error requirements of block LU factorization.

An important motivation for incorporating the Ozaki scheme into LU factorization comes from the High-Performance Linpack (HPL) benchmark. HPL solves a dense linear system using LU factorization with partial pivoting, with trailing matrix updates accounting for most of the arithmetic work for large problems. These updates are therefore a natural target for replacing conventional FP64 matrix multiplication by sequences of low-precision products. This possibility is already being explored in practice: NVIDIA's HPL implementation supports FP64 emulation using Ozaki schemes on Blackwell GPUs \cite{nvidiaHPL}. Such developments motivate a careful analysis of how the accuracy of the emulated matrix products affects the accuracy of the complete factorization.

The number of splits is central to this tradeoff. Using too many splits introduces additional matrix products and can diminish the performance advantage of low-precision hardware, whereas using too few can compromise the numerical reliability of the factorization. For an application such as HPL, the ultimate accuracy requirement concerns the computed solution, whose quality depends on the errors accumulated throughout the factorization and triangular solves. Consequently, imposing the same local accuracy requirement on every trailing matrix product may be unnecessarily restrictive. This motivates selecting the number of splits according to the contribution of each update to the overall factorization error, using information available as the computation progresses.


The remainder of the paper is organized as follows.
Section~\ref{sec:preliminary} reviews the Ozaki matrix multiplication algorithms proposed in~\cite{Ozaki2012} and~\cite{Ozaki2020}.
In Section~\ref{sec:algorithm-d}, we discuss how to choose $d$ for matrix multiplication and prove that the resulting choices meet the prescribed normwise and componentwise accuracy requirements.
Section~\ref{sec:algorithm-lu} then considers block LU factorization using Ozaki matrix multiplication for the trailing updates, derives corresponding choices of $d$, and proves the resulting accuracy guarantees.
Section~\ref{sec:experiments} presents numerical experiments illustrating the behavior of the proposed choices for block LU.

We use MATLAB-style notation to denote submatrices throughout the paper.
For example, $A_{i:j,k:l}$ denotes the submatrix of $A$ formed by rows $i,\ldots,j$ and columns $k,\ldots,l$.
We use $\LL$ and $\UU$ when the LU factors are viewed as block matrices or when referring to block submatrices, whereas $L$ and $U$ denote the same factors when individual entries, rows, or columns are considered.
Thus, the distinction between $\LL,\UU$ and $L,U$ is purely notational and reflects whether the argument is being carried out at the block or scalar level.
We use \(\norm{\cdot}\) for the \(2\)-norm.
For a matrix \(A=(A_{i,j})\), the notation \(\abs{A}\) denotes the componentwise absolute value of \(A\), i.e., \(\abs{A}=(\abs{A_{i,j}})\).

\section{The Ozaki matrix multiplication}
\label{sec:preliminary}
\subsection{Original Ozaki scheme}
The original Ozaki scheme~\cite{Ozaki2012} computes an accurate matrix product by decomposing the input matrices into a sequence of split matrices.
Let \(X\in\mathbb F^{m\times k}\) and \(Y\in\mathbb F^{k\times n}\).
The matrix \(X\) is split rowwise, whereas \(Y\) is split columnwise, i.e.,
\[
X=\sum_{r}X^{(r)},\qquad Y=\sum_{s}Y^{(s)},
\]
More specifically, starting with the residual \(R_X^{(1)}=X\), the \(r\)-th split of \(X\) is obtained by first determining, for each row \(i\), a power-of-two scaling factor from the largest entry of the current residual,
\[
\mu_i^{(r)}=\max_j\abs{(R_X^{(r)})_{ij}},\qquad
\sigma_i^{(r)}
=
2^{\left\lceil(\log_2 k-\log_2\macheps)/2\right\rceil}
2^{\left\lceil\log_2\mu_i^{(r)}\right\rceil},
\]
and then extracting its leading part by
\[
X^{(r)}
=
\flworking\!\left(
(R_X^{(r)}+\sigma^{(r)}e\trans)-\sigma^{(r)}e\trans
\right),
\qquad
R_X^{(r+1)}
=
\flworking\bigl(R_X^{(r)}-X^{(r)}\bigr),
\]
where \(e\) denotes the vector of all ones.
The same procedure is applied columnwise to \(Y\).
The scaling factor limits the number of significant bits contained in each split so that the corresponding split products can be evaluated without rounding error.
Repeating this procedure successively extracts the leading-to-trailing parts of the original matrices, while the splitting itself is error-free.
With sufficiently many splits, the residuals vanish and \(X\) and \(Y\) are represented exactly by their split matrices.
The resulting matrix product can then be transformed into a collection of error-free split products; together with error-free/accurate accumulation, the original scheme provides an error-free transformation of matrix multiplication.

\subsection{Fast Ozaki scheme}
For modern low-precision accelerators,~\cite{Ozaki2020}
adapted the splitting scheme to emulate high-precision matrix multiplication
using low-precision matrix multiplication units.
Consider $C=XY$, where $X\in\mathbb F^{m\times k}$ and $Y\in\mathbb F^{k\times n}$.
Let \(\macheps\), \(\machepsmedm\), and \(\machepslow\) denote the unit roundoffs of the working, medium, and low precisions, respectively.
In the FP64--FP32--FP16 implementation of~\cite{Ozaki2020}, these correspond to FP64, FP32, and FP16.
The input matrices and the final accumulation use \(\macheps\), the split matrices are stored in \(\machepslow\), and the products of the low-precision split matrices are accumulated in \(\machepsmedm\).
Define
\[
 \rho:=\left\lceil\log_2(\macheps^{-1})
 -\frac{\log_2(\machepsmedm^{-1})-\log_2 k}{2}\right\rceil,
 \qquad
 \beta:=\log_2(\macheps^{-1})-\rho
 =\left\lfloor\frac{\log_2(\machepsmedm^{-1})-\log_2 k}{2}\right\rfloor.
\]
Thus $\beta$ is the number of significant bits extracted at each splitting
step.  The choice of $\rho$, equivalently of $\beta$, ensures that the
inner products of the low-precision split vectors can be accumulated exactly
in the accumulation precision, provided the standard inner-product evaluation
is used.

The splitting procedure is summarized in Algorithm~\ref{alg:fast-ozaki-split}.
It is applied independently to every row of $X$ and every column of $Y$.
The scaling exponents are therefore row dependent for $X$ and column dependent
for $Y$.
In particular, if \(\bigl(c_X^{(r)}\bigr)_i\) and \(\bigl(c_Y^{(s)}\bigr)_j\) denote the exponents returned for row \(i\) and column \(j\), respectively, then
\begin{equation} \label{eq:splits-X-Y}
    X=\sum_r D_X^{(r)}X^{(r)},\qquad
    Y=\sum_s Y^{(s)}D_Y^{(s)},
\end{equation}
where
\[
    D_X^{(r)}=\diag(2^{\bigl(c_X^{(r)}\bigr)_1},\ldots,2^{\bigl(c_X^{(r)}\bigr)_m}),\qquad
    D_Y^{(s)}=\diag(2^{\bigl(c_Y^{(s)}\bigr)_1},\ldots,2^{\bigl(c_Y^{(s)}\bigr)_n}).
\]
In Section~\ref{sec:algorithm-d}, these power-of-two scaling matrices are absorbed into the split matrices, so \(X^{(r)}\) and \(Y^{(s)}\) there denote the splits in the original scale.

\begin{algorithm}[t]
\caption{Vector splitting in the fast Ozaki scheme~\cite{Ozaki2020}}
\label{alg:fast-ozaki-split}
\begin{algorithmic}[1]
\REQUIRE $x\in\mathbb F^k$, number of splits $d$, and $\rho$
\ENSURE low-precision splits $x^{(1)},\ldots,x^{(d)}$
        and scaling exponents $c^{(1)},\ldots,c^{(d)}$
\STATE $r\leftarrow x$
\FOR{$j=1,\ldots,d$}
    \STATE $\mu\leftarrow\max_i |r_i|$
    \STATE $c^{(j)}\leftarrow\lceil\log_2\mu\rceil$
    \STATE $\sigma\leftarrow 2^{\rho+c^{(j)}}$
    \FOR{$i=1,\ldots,k$}
        \STATE $x_{{\rm tmp},i}\leftarrow
        \flworking\bigl((r_i+\sigma)-\sigma\bigr)$
        \STATE $r_i\leftarrow\flworking(r_i-x_{{\rm tmp},i})$
        \STATE $x^{(j)}_i\leftarrow
        \fllow\bigl(2^{-c^{(j)}}x_{{\rm tmp},i}\bigr)$
    \ENDFOR
\ENDFOR
\end{algorithmic}
\end{algorithm}

After splitting, the fast matrix multiplication is summarized in
Algorithm~\ref{alg:fast-ozaki-gemm}.
For a prescribed $d$, only split products on the first $d$ anti-diagonals are retained, i.e., those satisfying $r+s\le d+1$.
Hence $d(d+1)/2$ low-precision GEMMs are required.
Each split product is evaluated with low-precision inputs and accumulation in the
medium precision with unit roundoff $\machepsmedm$; by the choice of $\rho$, these
split products are exact in the medium-precision accumulator.  The split
matrices themselves are stored in the low precision with unit roundoff
$\machepslow$.
The products are then rescaled by powers of two and accumulated in the working precision.
For standard working-precision accuracy, this final accumulation can be performed in the working precision, while NearSum can instead be used for correctly rounded accumulation~\cite{Ozaki2020}.

\begin{algorithm}[t]
\caption{Fast Ozaki matrix multiplication~\cite{Ozaki2020}}
\label{alg:fast-ozaki-gemm}
\begin{algorithmic}[1]
\REQUIRE $X\in\mathbb F^{m\times k}$, $Y\in\mathbb F^{k\times n}$,
         number of splits $d$
\ENSURE $\hat C\approx XY$
\STATE Compute $\rho$ (equivalently $\beta$) from $k$, $\macheps$, and
       $\machepsmedm$
\STATE Apply Algorithm~\ref{alg:fast-ozaki-split} rowwise to $X$ to obtain
       $\{X^{(r)},D_X^{(r)}\}_{r=1}^d$
\STATE Apply Algorithm~\ref{alg:fast-ozaki-split} columnwise to $Y$ to obtain
       $\{Y^{(s)},D_Y^{(s)}\}_{s=1}^d$
\STATE $\hat C\leftarrow 0$
\FOR{$r=1,\ldots,d$}
    \FOR{$s=1,\ldots,d+1-r$}
        \STATE $P^{(r,s)}\leftarrow
        X^{(r)}Y^{(s)}$
        \COMMENT{$\machepslow$ inputs with $\machepsmedm$ accumulation}
        \STATE $\hat C\leftarrow
        \flworking\bigl(\hat C+D_X^{(r)}P^{(r,s)}D_Y^{(s)}\bigr)$ \label{line:sum}
    \ENDFOR
\ENDFOR
\end{algorithmic}
\end{algorithm}

In the DP mode of~\cite{Ozaki2020}, \(d\) is selected before these matrix multiplications are performed.
Algorithm~3 of~\cite{Ozaki2020} successively generates the splits of \(X\) until
\[
    \flworking\!\left(
    (d+1)\abs{D_X^{(d)}X^{(d)}}(\abs{Y}e)
    \right)
    <
    \flworking\!\left(
    2\sqrt{k}\,\macheps\abs{X}(\abs{Y}e)
    \right)
\]
componentwise in the row index.
Algorithms~1 and 3 of~\cite{Ozaki2020} can therefore be performed simultaneously: the splits of \(X\) are generated until \(d\) is determined, after which the required splits of \(Y\) are formed.

The use of \(e\) avoids explicitly forming \(\abs{X}\abs{Y}\): the estimator requires only \(\abs{Y}e\) and matrix--vector products involving the splits of \(X\), giving quadratic cost for square matrices.
However, this compression discards columnwise information.
The stopping criterion naturally controls a rowwise quantity of the form
\[
    \abs{E}e\leq\bigO(\macheps)\abs{X}\abs{Y}e,
\]
rather than implying the componentwise bound
\[
    \abs{E}\leq\bigO(\macheps)\abs{X}\abs{Y}.
\]
Thus, the loss of a componentwise guarantee in the fast scheme is associated with the truncation and the rowwise selection of \(d\), not with the error-free splitting itself.
Moreover,~\cite[Algorithm~3]{Ozaki2020} selects \(d\) using a probabilistic error estimate, so the requested and achieved accuracies need not exactly coincide.

\begin{table}[t]
\centering
\caption{Computational costs of the fast Ozaki scheme for \(X\in\mathbb R^{m\times k}\) and \(Y\in\mathbb R^{k\times n}\).
Here \(c_{\rm mv}\) denotes the cost constant for a working-precision matrix--vector product, \(c_{\rm elem}\) the cost constant for elementwise operations, \(c_{\rm sum}\) the cost constant for rescaling and accumulating one split product in precision \(\macheps\), and \(c_{\rm gemm}\) the cost constant for a Tensor Core matrix multiplication with low precision inputs and medium precision accumulation.}
\label{tab:ozaki_cost}
\begin{tabular}{lcc}
\hline
&
\(\macheps\)
&
\(\machepsmedm/\machepslow\) on Tensor Cores\\
\hline
Error estimation
& \(2c_{\rm mv}k(n+m)\)
& \(0\)\\
Split \(X\) and select \(d\)
& \(c_{\rm elem}dmk+2c_{\rm mv}(d-1)mk\)
& \(0\)\\
Split \(Y\)
& \(c_{\rm elem}dkn\)
& \(0\)\\
Split products
& \(0\)
& \(c_{\rm gemm}\dfrac{d(d+1)}{2}mnk\)\\
Rescale and accumulate
& \(c_{\rm sum}\dfrac{d(d+1)}{2}mn\)
& \(0\)\\
\hline
\end{tabular}
\end{table}
\section{Strategies for selecting \(d\) in Ozaki multiplication}
\label{sec:algorithm-d}
As discussed in Section~\ref{sec:preliminary}, the split count \(d\) has a direct impact on performance.
We therefore explain how to choose \(d\) inexpensively under varying accuracy demands.
Assume that the computed result \(\hat{C}\) can be written as
\begin{equation} \label{eq:def-DXYgemm}
    \hat{C}=XY+\Delta(XY)^{(\mathtt{gemm})}(d,k),
\end{equation}
where $X\in\mathbb F^{m\times k}$ and $Y\in\mathbb F^{k\times n}$.
In the standard rounding error analysis for matrix--matrix multiplication, \(\hat{C}-XY\) typically depends on the inner dimension \(k\).
In this work we focus on Ozaki multiplication, for which \(\hat{C}-XY\) depends not only on \(k\) but also on \(d\).
We then consider two representative accuracy criteria:
\begin{enumerate}
\item Normwise stability: \(\abs{\Delta(XY)^{(\mathtt{gemm})}_{ij}(d,k)}\leq\varepsilon\norm{X_{i,:}}\norm{Y_{:,j}}\) for all \(i,j\);
\item Componentwise stability: \(\abs{\Delta(XY)^{(\mathtt{gemm})}(d,k)}\leq\varepsilon\abs{X}\abs{Y}\);
\end{enumerate}
where \(\varepsilon\) denotes a user-specified accuracy parameter.
Based on these two accuracy criteria, we propose corresponding strategies for choosing \(d\) to ensure the normwise and componentwise stability of matrix multiplication:
\begin{enumerate}
\item For normwise stability, we choose the smallest positive integer \(d\) such that
\begin{equation} \label{eq:d-normwise}
    4k(d+1)\theta^d\leq\frac{\varepsilon}{2}.
\end{equation}
Here \(\theta\) is determined by the inner dimension \(k\) and \(\machepsmedm\), as specified in Lemma~\ref{lem:bound-sum-X}.
\item For componentwise stability, we choose the smallest positive integer \(d\) such that
\begin{equation} \label{eq:d-componentwise}
    k(d+1)\theta^d\leq\frac{\varepsilon}{2}\xi,
\end{equation}
where \(\xi\) is any of the lower-bound estimates defined in Section~\ref{subsec:alo-component}.
\end{enumerate}
We then show that these choices guarantee the corresponding accuracy criteria for the fast Ozaki scheme.

Assuming that there is no overflow or problematic underflow, the downscaling and precision conversion in Algorithm~\ref{alg:fast-ozaki-split}, as well as the power-of-two rescaling in Line~\ref{line:sum} of Algorithm~\ref{alg:fast-ozaki-gemm}, are exact.
Therefore, the scaling matrices \(D_X^{(r)}\) and \(D_Y^{(s)}\) introduce no rounding errors and can be absorbed into the definitions of the split matrices.
Hence, without loss of generality, we take \(D_X^{(r)}=I\) and \(D_Y^{(s)}=I\) throughout the following error analysis.
With this convention, \(X^{(r)}\) and \(Y^{(s)}\) denote the splits in the original scale, and hence
\[
    X=\sum_r X^{(r)},\qquad
    Y=\sum_s Y^{(s)}.
\]
The fast Ozaki scheme can then be written as
\begin{equation} \label{eq:def-error-sum-trunc}
    C_d:=\sum_{r+s\leq d+1}X^{(r)}Y^{(s)},\qquad \Delta(XY)^{(\mathtt{sum})}(d):=\hat{C}-C_d,\qquad \Delta(XY)^{(\mathtt{trunc})}(d,k):=C_d-XY.
\end{equation}
Each \(X^{(r)}Y^{(s)}\) is computed exactly before the final working-precision accumulation.
Furthermore, \(\Delta(XY)^{(\mathtt{gemm})}(d,k)\) defined in~\eqref{eq:def-DXYgemm} can be written as
\begin{equation} \label{eq:DXYgemm-DXYsum+DXYtrunc}
    \Delta(XY)^{(\mathtt{gemm})}(d,k)=\Delta(XY)^{(\mathtt{sum})}(d)+\Delta(XY)^{(\mathtt{trunc})}(d,k).
\end{equation}

\subsection{Bound on \(\Delta(XY)^{(\mathtt{sum})}(d)\)}
We first consider \(\Delta(XY)^{(\mathtt{sum})}(d)\).
Since the split products are exact, this error only comes from their final accumulation in the working precision.
Therefore,
\begin{equation} \label{eq:DXY-sum}
\begin{split}
    \abs{\Delta(XY)^{(\mathtt{sum})}(d)}
    &\leq \gamma_{d(d+1)/2}\sum_{r+s\leq d+1}\abs{X^{(r)}Y^{(s)}}
    \leq \gamma_{d(d+1)/2}\sum_r\abs{X^{(r)}}\sum_s\abs{Y^{(s)}}.
\end{split}
\end{equation}
The following lemma bounds the absolute sums of the splits generated by the fast Ozaki scheme.

\begin{lemma} \label{lem:bound-sum-X}
Let \(X^{(r)}\) and \(Y^{(s)}\) denote the splits of \(X\) and \(Y\) generated by Algorithms~\ref{alg:fast-ozaki-split} and~\ref{alg:fast-ozaki-gemm}.
Let
\begin{equation*}
\begin{split}
    \rho&=\left\lceil\log_2(\macheps^{-1})-\frac{\log_2(\machepsmedm^{-1})-\log_2 k}{2}\right\rceil,\\
    \beta&:=\log_2(\macheps^{-1})-\rho
    =\left\lfloor\frac{\log_2(\machepsmedm^{-1})-\log_2 k}{2}\right\rfloor.
\end{split}
\end{equation*}
Define \(\theta:=2^{-\beta}\).
Assume that \(\rho\geq2\) and \(\beta\geq1\).
Also, assume that there is no overflow or problematic underflow.
Then
\begin{equation} \label{eq:sum-splits-X-Y}
    \sum_{r\geq1}\abs{X^{(r)}}\leq\kappa_{\rm sp}\abs{X},\qquad \sum_{s\geq1}\abs{Y^{(s)}}\leq\kappa_{\rm sp}\abs{Y},\qquad \kappa_{\rm sp}:=\frac{3+\theta}{1-\theta}.
\end{equation}
\end{lemma}

\begin{proof}
We only prove the result for \(X\), and then the argument for the columnwise splitting of \(Y\) is identical.
Considering row \(i\) at splitting step \(r\), by Algorithm~\ref{alg:fast-ozaki-split}, we have
\[
    \mu_i^{(r)}=\max_j\abs{(R_X^{(r)})_{ij}},\qquad
    \bigl(c_X^{(r)}\bigr)_i=\left\lceil\log_2\mu_i^{(r)}\right\rceil,\qquad
    \bigl(\sigma_X^{(r)}\bigr)_i =2^{\rho+\left(c_X^{(r)}\right)_i}.
\]
Since \(\abs{(R_X^{(r)})_{ij}}\leq2^{(c_X^{(r)})_i}\), we further have
\[
    \frac{\abs{(R_X^{(r)})_{ij}}}
    {(\sigma_X^{(r)})_i}
    \leq2^{-\rho}\leq\frac14.
\]
It follows that \(\flworking\!\left((R_X^{(r)})_{ij}+(\sigma_X^{(r)})_i\right)\) and \((\sigma_X^{(r)})_i\) are within a factor of two.
Together with Sterbenz's lemma, their subtraction in Algorithm~\ref{alg:fast-ozaki-split} is exact.

By \((\sigma_X^{(r)})_i=2^{\rho+(c_X^{(r)})_i}\), the spacing immediately below \((\sigma_X^{(r)})_i\) in the working precision is
\begin{equation*}
\begin{split}
    2^{\rho+(c_X^{(r)})_i-\log_2(\macheps^{-1})}
    &=2^{(c_X^{(r)})_i-\beta}
    =\theta 2^{(c_X^{(r)})_i},
\end{split}
\end{equation*}
whereas the spacing immediately above \((\sigma_X^{(r)})_i\) is twice this quantity.
Therefore, round-to-nearest gives
\begin{equation} \label{eq:residual-one-split}
    \abs{(R_X^{(r)})_{ij}-(X^{(r)})_{ij}}
    \leq \theta 2^{(c_X^{(r)})_i}.
\end{equation}
More precisely, the right-hand side can be replaced by \(\theta 2^{(c_X^{(r)})_i-1}\) when \((R_X^{(r)})_{ij}<0\).
Whenever \((X^{(r)})_{ij}\neq0\), the corresponding rounding intervals also show that \((R_X^{(r)})_{ij}\) and \((X^{(r)})_{ij}\) have the same sign and their magnitudes differ by at most a factor of two.
Hence Sterbenz's lemma shows that the residual update is also exact:
\[
    (R_X^{(r+1)})_{ij}
    =(R_X^{(r)})_{ij}-(X^{(r)})_{ij}.
\]
The case \((X^{(r)})_{ij}=0\) is trivial.
Thus both the extraction and the residual update in each splitting step are exact.

Taking the maximum over the row in~\eqref{eq:residual-one-split}, we obtain
\[
    \mu_i^{(r+1)}\leq\theta 2^{(c_X^{(r)})_i},
\]
which gives, by \(\theta=2^{-\beta}\),
\begin{equation} \label{eq:scale-decay-fast-ozaki}
    \bigl(c_X^{(r+1)}\bigr)_i
    \leq \bigl(c_X^{(r)}\bigr)_i-\beta.
\end{equation}

For a fixed nonzero entry \(X_{ij}\), let \(r_0\) be the first splitting step for which \((X^{(r_0)})_{ij}\neq0\).
Before this splitting step, \((R_X^{(r)})_{ij}=X_{ij}\).
From the rounding intervals at the \(r_0\) step,
\[
    2^{(c_X^{(r_0)})_i-\beta}\leq2\abs{X_{ij}},\qquad
    \abs{(R_X^{(r_0+1)})_{ij}}\leq\abs{X_{ij}}.
\]
Using~\eqref{eq:scale-decay-fast-ozaki}, for \(q\geq1\),
\[
    \abs{(R_X^{(r_0+q+1)})_{ij}}
    \leq2\theta^q\abs{X_{ij}}.
\]
Since each splitting step is exact, i.e.,
\(
    (X^{(r)})_{ij}
    =(R_X^{(r)})_{ij}-(R_X^{(r+1)})_{ij},
\)
we then obtain
\begin{equation*}
\begin{split}
    \sum_{r\geq1}\abs{(X^{(r)})_{ij}}
    &=\sum_{r\geq r_0}\abs{(X^{(r)})_{ij}} \\
    &\leq \abs{X_{ij}}
       +2\sum_{r\geq r_0+1}\abs{(R_X^{(r)})_{ij}} \\
    &\leq \abs{X_{ij}}+2\abs{X_{ij}}
       +4\sum_{q\geq1}\theta^q\abs{X_{ij}} \\
    &=\frac{3+\theta}{1-\theta}\abs{X_{ij}}.
\end{split}
\end{equation*}
This proves the first inequality in~\eqref{eq:sum-splits-X-Y} componentwise.
The second inequality follows analogously for \(Y\).
\end{proof}

Combining Lemma~\ref{lem:bound-sum-X} with~\eqref{eq:DXY-sum}, we immediately obtain a bound on \(\Delta(XY)^{(\mathtt{sum})}(d)\), which is stated in Lemma~\ref{lem:bound-XY-sum}.

\begin{lemma} \label{lem:bound-XY-sum}
Assume that \(\hat{C}\) is computed by the fast Ozaki scheme and \(\Delta(XY)^{(\mathtt{sum})}(d)\) is defined by~\eqref{eq:def-error-sum-trunc}.
Then
\begin{equation} \label{eq:bound-XY-sum}
    \abs{\Delta(XY)^{(\mathtt{sum})}(d)}\leq\gamma_{d(d+1)/2}\kappa_{\rm sp}^2\abs{X}\abs{Y},
\end{equation}
where \(\kappa_{\rm sp}\) is defined in~\eqref{eq:sum-splits-X-Y}.
\end{lemma}

\begin{remark} \label{remark:k-theta}
For the fast Ozaki scheme, \(\theta=2^{-\beta}\), where
\[
    \beta=\left\lfloor\frac{\log_2(\machepsmedm^{-1})-\log_2 k}{2}\right\rfloor.
\]
Thus, for the inner dimension \(k\), \(\theta\) is typically much smaller than one.
Consequently,
\[
    \kappa_{\rm sp}=\frac{3+\theta}{1-\theta}=\bigO(1).
\]
For example, when \(\machepsmedm=2^{-24}\), corresponding to FP32 accumulation, some representative values are
\[
\begin{array}{c|cccccccc}
    k & 16 & 32 & 64 & 128 & 256 & 512 & 1024 & 2^{21} \\
    \hline
    \theta & 2^{-10} & 2^{-9} & 2^{-9} & 2^{-8} & 2^{-8} & 2^{-7} & 2^{-7} & 2^{-1}
\end{array}
\]
In particular, \(\theta=1/2\) first occurs among powers of two at \(k=2^{21}\), which is already a very large inner dimension.
This further means that \(\kappa_{\rm sp}\leq 7\) holds for all \(k\leq 2^{21}\).
\end{remark}

\subsection{Normwise stability for the strategy~\eqref{eq:d-normwise}}
\label{subsec:alo-norm}
We next bound the truncation error caused by omitting the split products with \(r+s>d+1\).
The following lemma gives bounds for both the truncation error \(\Delta(XY)^{(\mathtt{trunc})}(d,k)\) and the total error \(\Delta(XY)^{(\mathtt{gemm})}(d,k)\).

\begin{theorem} \label{lem:bound-XY-normwise}
Let \(\hat C\) be computed by the fast Ozaki scheme with a given splits number \(d\), and let \(\Delta(XY)^{(\mathtt{sum})}(d)\) and \(\Delta(XY)^{(\mathtt{trunc})}(d,k)\) be defined by~\eqref{eq:def-error-sum-trunc}.
Then, for all \(i\) and \(j\),
\begin{equation} \label{eq:bound-XY-trunc-entry}
    \abs{\Delta(XY)_{ij}^{(\mathtt{trunc})}(d,k)}\leq4k(d+1)\theta^d\norm{X_{i,:}}\norm{Y_{:,j}},
\end{equation}
and
\begin{equation} \label{eq:bound-XY-gemm-entry}
    \abs{\Delta(XY)^{(\mathtt{gemm})}_{ij}(d,k)}\leq\left(4k(d+1)\theta^d+\gamma_{d(d+1)/2}\kappa_{\rm sp}^2\right)\norm{X_{i,:}}\norm{Y_{:,j}}.
\end{equation}
\end{theorem}

\begin{proof}
For each row of \(X\) and each column of \(Y\), define
\begin{equation}
    \widetilde{x}_i:=2^{\lceil\log_2\max_j|X_{ij}|\rceil},\qquad \widetilde{y}_j:=2^{\lceil\log_2\max_i|Y_{ij}|\rceil}.
\end{equation}
By the residual bound and the scale decay in the proof of Lemma~\ref{lem:bound-sum-X}, the successive splits and residuals satisfy
\begin{equation} \label{eq:split-envelope-X-Y}
    \abs{(X^{(r)})_{ij}}\leq\theta^{r-1}\widetilde{x}_i,\qquad \abs{(Y^{(s)})_{ij}}\leq\theta^{s-1}\widetilde{y}_j,
\end{equation}
and
\begin{equation} \label{eq:residual-envelope-X-Y}
    \abs{(R_X^{(r+1)})_{ij}}\leq\theta^r\widetilde{x}_i,\qquad \abs{(R_Y^{(s+1)})_{ij}}\leq\theta^s\widetilde{y}_j.
\end{equation}
Since \(X=\sum_{r\geq1}X^{(r)}\) and \(Y=\sum_{s\geq1}Y^{(s)}\), the omitted part can be written as
\begin{equation} \label{eq:trunc-decomposition}
    XY-C_d=\sum_{r=1}^dX^{(r)}R_Y^{(d+2-r)}+R_X^{(d+1)}Y.
\end{equation}
Using~\eqref{eq:split-envelope-X-Y} and~\eqref{eq:residual-envelope-X-Y}, we obtain
\begin{equation} \label{eq:bound-XY-trunc-component}
\begin{split}
    \abs{\Delta(XY)_{ij}^{(\mathtt{trunc})}(d,k)}
    &=\abs{(C_d-XY)_{ij}}\\
    &\leq k(d+1)\theta^d\widetilde{x}_i\widetilde{y}_j\\
    &<4k(d+1)\theta^d\norm{X_{i,:}}\norm{Y_{:,j}},
\end{split}
\end{equation}
where the last inequality follows from \(\widetilde{x}_i<2\norm{X_{i,:}}\) and \(\widetilde{y}_j<2\norm{Y_{:,j}}\).
This inequality~\eqref{eq:bound-XY-trunc-component} proves~\eqref{eq:bound-XY-trunc-entry}.
Combining Lemma~\ref{lem:bound-XY-sum} and~\eqref{eq:bound-XY-trunc-entry} with~\eqref{eq:DXYgemm-DXYsum+DXYtrunc},
we have~\eqref{eq:bound-XY-gemm-entry}.
\end{proof}

\begin{remark}
Although \(\gamma_{d(d+1)/2}\) involved in~\eqref{eq:bound-XY-gemm-entry} increases with \(d\), the split count \(d\) is typically small in the fast Ozaki scheme.
Moreover, as discussed in Remark~\ref{remark:k-theta}, \(\kappa_{\rm sp}=\bigO(1)\) for practical inner dimensions \(k\).
It is therefore reasonable in the regime of interest to assume that
\begin{equation} \label{eq:sum-half-eps}
    \gamma_{d(d+1)/2}\kappa_{\rm sp}^2\leq\frac{\varepsilon}{2}.
\end{equation}
Under~\eqref{eq:sum-half-eps}, it remains to control the truncation error by requiring
\begin{equation}
    4k(d+1)\theta^d\leq\frac{\varepsilon}{2}.
\end{equation}
Hence, the choice~\eqref{eq:d-normwise} and Lemma~\ref{lem:bound-XY-normwise} can guarantee
\begin{equation}
    \abs{\Delta(XY)^{(\mathtt{gemm})}_{ij}(d,k)}\leq\varepsilon\norm{X_{i,:}}\norm{Y_{:,j}}
\end{equation}
for all \(i\) and \(j\).
\end{remark}

\subsection{Componentwise stability for the strategy~\eqref{eq:d-componentwise}}
\label{subsec:alo-component}
We finally consider componentwise stability.
The main difficulty is that the truncation bound in~\eqref{eq:bound-XY-trunc-component} is expressed in terms of the row and column norms, whereas the desired componentwise scale is \(\abs{X}\abs{Y}\).
We therefore construct lower-bound estimates for \(\abs{X}\abs{Y}\) that can be evaluated much more cheaply than forming \(\abs{X}\abs{Y}\).
Recall that \(\widetilde{x}_i\) and \(\widetilde{y}_j\) are defined in~\eqref{eq:split-envelope-X-Y}.
For nonzero \(\widetilde{x}_i\) and \(\widetilde{y}_j\), define
\begin{equation} \label{eq:normalized-z-w}
    z_{il}:=\frac{\abs{X_{il}}}{\widetilde{x}_i},\qquad w_{l j}:=\frac{\abs{Y_{l j}}}{\widetilde{y}_j}.
\end{equation}
Then \(0\leq z_{il},w_{l j}\leq1\).
We consider the following three choices of the lower-bound estimate \(\xi\).
\begin{enumerate}
\item Define
\begin{equation}
    m_l^{(X)}:=\min_{i:\widetilde{x}_i>0}z_{il},\qquad m_l^{(Y)}:=\min_{j:\widetilde{y}_j>0}w_{l j},\qquad
    \xi_{\rm col}:=\sum_{l=1}^k m_l^{(X)}m_l^{(Y)},
\end{equation}
and
\begin{equation}
    s_X:=\min_{i:\widetilde{x}_i>0}\sum_{l=1}^k z_{il},\qquad s_Y:=\min_{j:\widetilde{y}_j>0}\sum_{l=1}^k w_{l j}.
\end{equation}
Then take
\begin{equation} \label{eq:xi-componentwise}
    \xi^{(1)}:=\max\{\xi_{\rm col},\xi_{\rm sum}\},\qquad \xi_{\rm sum}:=\max\{0,s_X+s_Y-k\}.
\end{equation}
This choice requires \(\bigO(mk+kn)\) work, while \(\widetilde{x}_i\) and \(\widetilde{y}_j\) are also required by the splitting process and can therefore be reused.

\item For each row of \(X\) and each column of \(Y\), let
\begin{equation}
    p_i\in\arg\max_{1\leq l\leq k}z_{il},\qquad q_j\in\arg\max_{1\leq l\leq k}w_{l j}.
\end{equation}
Define
\begin{equation} \label{eq:xi-one-index}
    \xi^{(2)}:=\min_{i,j}\sum_{l\in\{p_i,q_j\}}z_{il}w_{l j},
\end{equation}
where a repeated index is included only once.
The indices \(p_i\) and \(q_j\) can be obtained when \(\widetilde{x}_i\) and \(\widetilde{y}_j\) are computed, and evaluating~\eqref{eq:xi-one-index} therefore requires only \(\bigO(mn)\) additional work.

\item Let \(S_i\) contain the indices of the \(s\) largest entries in \(\{z_{il}\}_{l=1}^k\), and let \(T_j\) contain the indices of the \(s\) largest entries in \(\{w_{l j}\}_{l=1}^k\).
Define
\begin{equation} \label{eq:xi-top-s}
    \xi^{(3)}:=\min_{i,j}\sum_{l\in S_i\cup T_j}z_{il}w_{l j}.
\end{equation}
For a fixed small \(s\), evaluating~\eqref{eq:xi-top-s} after the index sets have been obtained requires \(\bigO(smn)\) additional work.
It should be noted that the second choice is the special case \(s=1\) of this choice.
\end{enumerate}

We then demonstrate that each of these three choices can ensure componentwise stability.

\begin{theorem} \label{lem:bound-XY-componentwise}
Let \(\hat C\) be computed by the fast Ozaki scheme with a given split number \(d\), and let \(\Delta(XY)^{(\mathtt{sum})}(d)\) and \(\Delta(XY)^{(\mathtt{trunc})}(d,k)\) be defined by~\eqref{eq:def-error-sum-trunc}.
Let \(\xi\) be any one of \(\xi^{(1)}\), \(\xi^{(2)}\), and \(\xi^{(3)}\) defined by~\eqref{eq:xi-componentwise}, \eqref{eq:xi-one-index}, and~\eqref{eq:xi-top-s}, respectively.
Then
\begin{equation} \label{eq:lower-reference-componentwise}
    \abs{X}\abs{Y}\geq\xi\,\widetilde{x}\widetilde{y}^{T},
\end{equation}
where \(\widetilde{x}=(\widetilde{x}_i)\) and \(\widetilde{y}=(\widetilde{y}_j)\) are defined in~\eqref{eq:split-envelope-X-Y}.
If \(\xi>0\), then
\begin{equation} \label{eq:bound-XY-trunc-cw}
    \abs{\Delta(XY)^{(\mathtt{trunc})}(d,k)}\leq\frac{k(d+1)\theta^d}{\xi}\abs{X}\abs{Y},
\end{equation}
and
\begin{equation} \label{eq:bound-XY-gemm-cw}
    \abs{\Delta(XY)^{(\mathtt{gemm})}(d,k)}\leq\left(\frac{k(d+1)\theta^d}{\xi}+\gamma_{d(d+1)/2}\kappa_{\rm sp}^2\right)\abs{X}\abs{Y}.
\end{equation}
\end{theorem}

\begin{proof}
For every \(i\), \(j\) with \(\widetilde{x}_i\widetilde{y}_j>0\), we have
\begin{equation} \label{eq:normalized-XY-componentwise}
    (\abs{X}\abs{Y})_{ij}=\widetilde{x}_i\widetilde{y}_j\sum_{l=1}^k z_{il}w_{l j}.
\end{equation}
We first verify that each of the three choices is a lower bound for the normalized inner product in~\eqref{eq:normalized-XY-componentwise}.
For \(\xi^{(1)}\), the definitions of \(m_l^{(X)}\) and \(m_l^{(Y)}\) give
\begin{equation}
    \sum_{l=1}^k z_{il}w_{l j}\geq\sum_{l=1}^k m_l^{(X)}m_l^{(Y)}=\xi_{\rm col}.
\end{equation}
Moreover, from \(0\leq z_{il}\), \(w_{l j}\leq1\), we have \(z_{il}w_{l j}\geq z_{il}+w_{l j}-1\), and hence
\begin{equation}
\begin{split}
    \sum_{l=1}^k z_{il}w_{l j}
    &\geq\sum_{l=1}^k z_{il}+\sum_{l=1}^k w_{l j}-k
    \geq s_X+s_Y-k.
\end{split}
\end{equation}
Since the left-hand side is nonnegative, we then have
\begin{equation}
    \sum_{l=1}^k z_{il}w_{l j}\geq\xi_{\rm sum},
\end{equation}
and therefore \(\sum_{l=1}^k z_{il}w_{l j}\geq\xi^{(1)}\).
For \(\xi^{(2)}\), all terms are nonnegative, thus
\begin{equation}
    \sum_{l=1}^k z_{il}w_{l j}\geq\sum_{l\in\{p_i,q_j\}}z_{il}w_{l j}\geq\xi^{(2)}.
\end{equation}
Similarly, for \(\xi^{(3)}\),
\begin{equation}
    \sum_{l=1}^k z_{il}w_{l j}\geq\sum_{l\in S_i\cup T_j}z_{il}w_{l j}\geq\xi^{(3)}.
\end{equation}
Thus~\eqref{eq:lower-reference-componentwise} holds for any of the three choices of \(\xi\).
If \(\widetilde{x}_i\widetilde{y}_j=0\), both sides of~\eqref{eq:lower-reference-componentwise} are zero, so the inequality remains valid.

From~\eqref{eq:bound-XY-trunc-component}, we obtain
\begin{equation}
    \abs{\Delta(XY)^{(\mathtt{trunc})}(d,k)}\leq k(d+1)\theta^d\widetilde{x}\widetilde{y}^{T}.
\end{equation}
When \(\xi>0\), combining this inequality with~\eqref{eq:lower-reference-componentwise} yields~\eqref{eq:bound-XY-trunc-cw}.
Finally, Lemma~\ref{lem:bound-XY-sum} and~\eqref{eq:DXYgemm-DXYsum+DXYtrunc} give~\eqref{eq:bound-XY-gemm-cw}.
\end{proof}

\begin{remark}
All three choices of \(\xi\) avoid forming the full matrix product \(\abs{X}\abs{Y}\).
The first choice only costs \(\bigO(mk+kn)\), whereas the second and third choices require \(\bigO(mn)\) and \(\bigO(smn)\) additional work, respectively, once the required indices are available.
As in the normwise case, we only need to control the truncation error by requiring
\begin{equation}
    \frac{k(d+1)\theta^d}{\xi}\leq\frac{\varepsilon}{2},
\end{equation}
which is exactly the selection criterion~\eqref{eq:d-componentwise}.
Therefore, if \(\xi>0\), Theorem~\ref{lem:bound-XY-componentwise} guarantees
\begin{equation}
    \abs{\Delta(XY)^{(\mathtt{gemm})}(d,k)}\leq\varepsilon\abs{X}\abs{Y}.
\end{equation}
If \(\xi=0\), the corresponding lower reference is inconclusive and~\eqref{eq:d-componentwise} cannot certify a finite \(d\), even though the actual componentwise error may still satisfy the desired bound.
\end{remark}

\section{Strategies for selecting \(d\) in block LU}
\label{sec:algorithm-lu}
In this section, we consider using the fast Ozaki scheme in LU factorization.
We begin by revisiting the block LU factorization algorithm summarized in Algorithm~\ref{alg:block-LU}.
Within this algorithm, only Line~\ref{line:update_trail} involves a matrix--matrix multiplication.

\begin{algorithm}[!tb]
\caption{The block LU factorization}
\label{alg:block-LU}
\begin{algorithmic}[1]
\REQUIRE A block size \(n_b\), and a given matrix \(\Acal\in\mathbb F^{m\times n}\) 
with \(m\geq n\),  \(m = qn_b\), and \(n = pn_b\).
\ENSURE A lower trapezoidal matrix \(\LL\in\mathbb F^{m\times n}\)
and an upper triangular matrix \(\UU\in\mathbb F^{n\times n}\)
satisfying \(\Acal = \LL\UU\).

\FOR{\(i = 1:p-1\)}
    \STATE Compute the LU factorization \(\Acal_{i:q, i} = \LL_{i:q, i}
    \UU_{i,i}\). \label{line:lu:sub-lu}
    \STATE Update \(\UU_{i, i+1:p}\) by solving the lower triangular 
    system \(\LL_{i,i} \UU_{i, i+1:p} = \Acal_{i, i+1:p}\). \label{line:lu:tr-solver}
    \STATE Update trailing submatrix by \(\Acal_{i+1:q, i+1:p} = 
    \Acal_{i+1:q, i+1:p} - \LL_{i+1:q, i} \UU_{i, i+1:p}\). \label{line:update_trail}
\ENDFOR
\STATE Compute the LU factorization \(\Acal_{p:q, p} = \LL_{p:q, p}
\UU_{p,p}\).
\end{algorithmic}
\end{algorithm}

For each iteration in Algorithm~\ref{alg:block-LU}, we assume that the LU factorization and the triangular solve in Lines~\ref{line:lu:sub-lu} and~\ref{line:lu:tr-solver} are both performed in \(\macheps\).
Also, we assume that the matrix--matrix subtraction in Line~\ref{line:update_trail} is performed in \(\macheps\).
Thus, from~\cite[Theorem 2]{CMS2026}, we have the following corollary.

\begin{corollary} \label{cor:block-lu}
    Let \(\LL\in\mathbb{F}^{m\times n}\) and \(\UU\in\mathbb{F}^{n\times n}\)
    denote a lower trapezoidal matrix and an upper triangular matrix, respectively.
    Assume that \(\Acal = \LL\UU\) is computed by
    Algorithm~\ref{alg:block-LU}.
    Also, assume that Line~\ref{line:update_trail} in
    Algorithm~\ref{alg:block-LU} satisfy
    \begin{align} \label{eq:cor:blocklu:update-trail}
        \Acalhat^{(i-1)}_{i+1:q,i+1:p}
        - \LLhat_{i+1:q,i}\UUhat_{i,i+1:p}
        &= \Acalhat^{(i)}_{i+1:q,i+1:p}
        + \bigl(\Delta\bigl(\LL_{i+1:q,i}\UU_{i,i+1:p}\bigr)^{(\mathtt{gemm})}(n_b)\bigr)^{(i)}
        \notag\\
        &\quad
        + \bigl(\Delta\Acal^{(i-1)}_{i+1:q,i+1:p}\bigr)^{(\mathtt{gemm})},
    \end{align}
    where \(\Acalhat^{(0)}_{i+1:q,i+1:p}:=\Acal_{i+1:q,i+1:p}\) and \(\abs{\bigl(\Delta\Acal^{(i-1)}_{i+1:q,i+1:p}\bigr)^{(\mathtt{gemm})}}\leq\macheps\abs{\Acalhat^{(i-1)}_{i+1:q,i+1:p}}\).
    Then there exists \(\Delta(\LL\UU)^{(\mathtt{getrf})}(n,n_b)\) such that
    \begin{equation} \label{eq:blocklu:getrf}
        \Acal + \Delta(\LL\UU)^{(\mathtt{getrf})}(n,n_b) = \LLhat\UUhat,
    \end{equation}
    where, for \(i=1,\ldots,q\), \(j=1,\ldots,p\), and \(s=\min\{i,j\}\),
    \begin{equation} \label{eq:blocklu:getrf-trackable}
    \begin{split}
        \abs{\Delta(\LL\UU)^{(\mathtt{getrf})}_{i,j}(n,n_b)}
        \leq{}&(1+\gamma_{s-1})\left(
        \sum_{k=1}^{s-1}
        \abs{\bigl(\Delta(\LL_{i,k}\UU_{k,j})^{(\mathtt{gemm})}(n_b)\bigr)^{(k)}}
        +\gamma_{n_b}\abs{\LLhat_{i,s}}\abs{\UUhat_{s,j}}\right)
        \\
        &+\gamma_{s-1}
        \abs{\LLhat_{i,1:s}}\abs{\UUhat_{1:s,j}}.
    \end{split}
    \end{equation}
\end{corollary}

\begin{proof}
From~\cite[Theorem 2]{CMS2026}, before bounding the subtraction errors, the perturbation in the \(i\)th block column satisfies
\begin{equation} \label{eq:thm-proof:firstbound}
\begin{split}
    \abs{\Delta(\LL\UU)^{(\mathtt{getrf})}_{i:q,i}(n,n_b)}
    \leq{}&\sum_{k=1}^{i-1}
    \abs{\bigl(\Delta(\LL_{i:q,k}\UU_{k,i})^{(\mathtt{gemm})}(n_b)\bigr)^{(k)}}
    \\
    &+\sum_{k=0}^{i-2}
    \abs{\bigl(\Delta\Acal^{(k)}_{i:q,i}\bigr)^{(\mathtt{gemm})}}
    +\gamma_{n_b}\abs{\LLhat_{i:q,i}}\abs{\UUhat_{i,i}}.
\end{split}
\end{equation}
We first bound the second sum.
For \(k=1,\ldots,i-1\), restricting~\eqref{eq:cor:blocklu:update-trail} to the \(i\)th block column and using
\(\abs{(\Delta\Acal^{(k-1)})^{(\mathtt{gemm})}}\leq\macheps\abs{\Acalhat^{(k-1)}}\) gives
\begin{equation} \label{eq:blocklu-A-recurrence}
\begin{split}
    (1-\macheps)\abs{\Acalhat^{(k-1)}_{i:q,i}}
    \leq{}&\abs{\Acalhat^{(k)}_{i:q,i}}
    +\abs{\LLhat_{i:q,k}}\abs{\UUhat_{k,i}}
    \\
    &+\abs{\bigl(\Delta(\LL_{i:q,k}\UU_{k,i})^{(\mathtt{gemm})}(n_b)\bigr)^{(k)}}.
\end{split}
\end{equation}
After the first \(i-1\) trailing updates, \(\Acalhat^{(i-1)}_{i:q,i}\) is the input to the panel factorization in the \(i\)th block column.
The backward error of this panel factorization therefore gives
\begin{equation} \label{eq:blocklu-terminal-A}
    \abs{\Acalhat^{(i-1)}_{i:q,i}}
    \leq(1+\gamma_{n_b})\abs{\LLhat_{i:q,i}}\abs{\UUhat_{i,i}}.
\end{equation}
Applying~\eqref{eq:blocklu-A-recurrence} backwards from~\eqref{eq:blocklu-terminal-A}, and using
\((1-\macheps)^{-j}\leq 1+\gamma_j\), yields
\begin{equation} \label{eq:blocklu-subtraction-errors}
\begin{split}
    \sum_{k=0}^{i-2}
    \abs{\bigl(\Delta\Acal^{(k)}_{i:q,i}\bigr)^{(\mathtt{gemm})}}
    \leq{}&\gamma_{i-1}\Bigg(
    \sum_{k=1}^{i-1}
    \abs{\bigl(\Delta(\LL_{i:q,k}\UU_{k,i})^{(\mathtt{gemm})}(n_b)\bigr)^{(k)}}
    \\
    &\qquad+\abs{\LLhat_{i:q,1:i-1}}\abs{\UUhat_{1:i-1,i}}
    +(1+\gamma_{n_b})\abs{\LLhat_{i:q,i}}\abs{\UUhat_{i,i}}
    \Bigg).
\end{split}
\end{equation}
Substituting~\eqref{eq:blocklu-subtraction-errors} into~\eqref{eq:thm-proof:firstbound} and using
\[
    \abs{\LLhat_{i:q,1:i-1}}\abs{\UUhat_{1:i-1,i}}
    +\abs{\LLhat_{i:q,i}}\abs{\UUhat_{i,i}}
    =\abs{\LLhat_{i:q,1:i}}\abs{\UUhat_{1:i,i}}
\]
gives~\eqref{eq:blocklu:getrf-trackable} with \(s=j\) for blocks on and below the diagonal.
For a block above the diagonal, the corresponding identity in~\cite[Theorem~2]{CMS2026} has the same structure, with the panel-factorization error replaced by the triangular-solve error.
Since both the panel-factorization and triangular-solve errors are bounded by \(\gamma_{n_b}\abs{\LLhat_{i,s}}\abs{\UUhat_{s,j}}\), the same argument as \(s=i\) also proves~\eqref{eq:blocklu:getrf-trackable} for these blocks.
\end{proof}

To apply Corollary~\ref{cor:block-lu} to Algorithm~\ref{alg:block-LU} with the fast Ozaki scheme,
we first specify how the errors in the trailing-matrix update are assigned
to the two error terms in \eqref{eq:cor:blocklu:update-trail}.
For simplicity, let
\[
    \hat{C}^{(i)} = \flworking\bigl(\LLhat_{i+1:q,i}\UUhat_{i,i+1:p}\bigr)
\]
denote the matrix product returned by the matrix multiplication kernel.
The subsequent matrix subtraction is performed componentwise, and hence
\begin{equation*}
\begin{split}
    \abs{\flworking\bigl(\Acalhat^{(i-1)}_{i+1:q,i+1:p}-\hat{C}^{(i)}\bigr)-\bigl(\Acalhat^{(i-1)}_{i+1:q,i+1:p}-\hat{C}^{(i)}\bigr)}
    \leq{}& \macheps\abs{\Acalhat^{(i-1)}_{i+1:q,i+1:p}}+\macheps\abs{\hat{C}^{(i)}}.
\end{split}
\end{equation*}
Therefore, the contribution of the subtraction associated with
\(\Acalhat^{(i-1)}_{i+1:q,i+1:p}\) is included in the term \(\bigl(\Delta\Acal^{(i-1)}_{i+1:q,i+1:p}\bigr)^{(\mathtt{gemm})}\) satisfying
\[
    \abs{\bigl(\Delta\Acal^{(i-1)}_{i+1:q,i+1:p}\bigr)^{(\mathtt{gemm})}}
    \leq\macheps\abs{\Acalhat^{(i-1)}_{i+1:q,i+1:p}}.
\]
The contribution of the subtraction associated with the computed matrix
product \(\hat{C}^{(i)}\), on the other hand, is included in
\[
    \bigl(
        \Delta(
        \LL_{i+1:q,i}\UU_{i,i+1:p}
        )^{(\mathtt{gemm})}(n_b)
    \bigr)^{(i)}.
\]

Suppose now that the product \(\LLhat_{i+1:q,i}\UUhat_{i,i+1:p}\) is evaluated by the fast Ozaki scheme.
Since the split matrix products are computed exactly, the error in \(\hat{C}^{(i)}\) consists of the summation error and the truncation error.
Thus, we write
\[
    \hat{C}^{(i)}=\LLhat_{i+1:q,i}\UUhat_{i,i+1:p}
    +\bigl(\Delta(\LL_{i+1:q,i}\UU_{i,i+1:p})^{(\mathtt{sum})}(d)\bigr)^{(i)}
    +\bigl(\Delta(\LL_{i+1:q,i}\UU_{i,i+1:p})^{(\mathtt{trunc})}(d,n_b)\bigr)^{(i)}.
\]
Consequently, the matrix-product error appearing in
\eqref{eq:cor:blocklu:update-trail} is given explicitly by
\begin{equation} \label{eq:bound-DLU}
\begin{split}
    &\abs{\bigl(\Delta(\LL_{i+1:q,i}\UU_{i,i+1:p})^{(\mathtt{gemm})}(n_b)\bigr)^{(i)}} \\
    &\qquad\leq\abs{\bigl(\Delta(\LL_{i+1:q,i}\UU_{i,i+1:p})^{(\mathtt{sum})}(d)\bigr)^{(i)}}
    +\abs{\bigl(\Delta(\LL_{i+1:q,i}\UU_{i,i+1:p})^{(\mathtt{trunc})}(d,n_b)\bigr)^{(i)}}
    +\macheps\abs{\hat{C}^{(i)}} \\
    &\qquad\leq(1+\macheps)\abs{\bigl(\Delta(\LL_{i+1:q,i}\UU_{i,i+1:p})^{(\mathtt{sum})}(d)\bigr)^{(i)}}
    +(1+\macheps)\abs{\bigl(\Delta(\LL_{i+1:q,i}\UU_{i,i+1:p})^{(\mathtt{trunc})}(d,n_b)\bigr)^{(i)}} \\
    &\qquad\quad+\macheps\abs{\LLhat_{i+1:q,i}\UUhat_{i,i+1:p}}.
\end{split}
\end{equation}
By Lemma~\ref{lem:bound-XY-sum}, the summation error in each trailing update has already been bounded by
\begin{equation} \label{eq:lu-sum-bound}
    \abs{\bigl(\Delta(\LL_{i+1:q,i}\UU_{i,i+1:p})^{(\mathtt{sum})}(d)\bigr)^{(i)}}
    \leq\gamma_{d_i(d_i+1)/2}\kappa_{\rm sp}^2\abs{\LLhat_{i+1:q,i}}\abs{\UUhat_{i,i+1:p}}.
\end{equation}
Hence, the only new error that needs to be controlled through the choice of \(d_i\) is
\[
    \abs{\bigl(\Delta(\LL_{i+1:q,i}\UU_{i,i+1:p})^{(\mathtt{trunc})}(d,n_b)\bigr)^{(i)}}.
\]
We discuss normwise and componentwise choices separately.

\subsection{Normwise stability}
\label{subsec:lu-normwise}
In this subsection, \(L\) and \(\LL\) denote the same matrix, and similarly for \(U\) and \(\UU\).
We use \(\LL\) and \(\UU\) for blocks and block submatrices, and \(L\) and \(U\) for scalar entries and scalar row or column segments.
Thus, both notations are used in the proof below: block-level relations are written in terms of \(\LL\) and \(\UU\), whereas entrywise relations are written in terms of \(L\) and \(U\).
Similarly, \(\Delta(\LL\UU)^{(\mathtt{getrf})}\) denotes the backward-error matrix at the block level, while \(\Delta(LU)^{(\mathtt{getrf})}_{rs}\) denotes its scalar \((r,s)\) entry.

\subsubsection{A fixed selection with essentially no additional cost}
\label{subsubsec:lu-normwise-fixed}

A direct strategy is to use the normwise selection rule for matrix multiplication without exploiting any additional information from the LU factorization.
Indeed, the truncation error in each trailing update is exactly of the form analyzed in Theorem~\ref{lem:bound-XY-normwise} with inner dimension \(n_b\).
Let \(d\) be the smallest positive integer satisfying
\begin{equation} \label{eq:lu-d-normwise}
    4n_b(d+1)\theta^d\leq\frac{\varepsilon}{2}.
\end{equation}
Here both \(n_b\) and \(\varepsilon\) are fixed throughout the LU factorization, and \(\theta\) depends only on the fixed inner dimension \(n_b\).
Therefore, the same \(d\) is used for all trailing updates and needs to be determined only once before the factorization.
The selection consists only of checking the scalar inequality~\eqref{eq:lu-d-normwise} and hence introduces essentially no additional computational cost.
Then, for the \(i\)th trailing update, we have
\begin{equation} \label{eq:lu-trunc-normwise}
    \abs{\bigl(\Delta(\LL_{i+1:q,i}\UU_{i,i+1:p})_{rs}^{(\mathtt{trunc})}(d,n_b)\bigr)^{(i)}}
    \leq\frac{\varepsilon}{2}\norm{(\LLhat_{i+1:q,i})_{r,:}}\norm{(\UUhat_{i,i+1:p})_{:,s}}.
\end{equation}
Together with~\eqref{eq:lu-sum-bound}, this gives the same normwise accuracy for the matrix multiplication used in each trailing update as in Section~\ref{subsec:alo-norm}.
The following result shows that these local bounds also have the correct scale for the final LU factorization.
Although the adaptive result below is closely related, the fixed strategy admits a sharper bound because its local errors scale with the norms of the current block segments and can therefore be bounded by the Cauchy--Schwarz inequality without an iteration-count factor.

\begin{theorem} \label{thm:lu-normwise-d}
Assume that \(d\) is chosen by~\eqref{eq:lu-d-normwise} and that \(\varepsilon\) satisfies~\eqref{eq:sum-half-eps} for this \(d\).
Then there exists \(\Delta(LU)^{(\mathtt{getrf})}(n,n_b)=\Delta(\LL\UU)^{(\mathtt{getrf})}(n,n_b)\) such that~\eqref{eq:blocklu:getrf} holds with, for any scalar entry \((r,t)\) in block \((i,j)\), and \(s=\min\{i,j\}\),
\begin{equation} \label{eq:lu-gemm-normwise-final}
\begin{split}
    \abs{\Delta(LU)^{(\mathtt{getrf})}_{rt}(n,n_b)}
    \leq{}&2(1+\gamma_{s-1})
    \bigl(\varepsilon+\gamma_{n_b}+\gamma_{s-1}\bigr)
    \norm{\hat L_{r,:}}\norm{\hat U_{:,t}}.
\end{split}
\end{equation}
In particular, if \(\varepsilon=\bigO(\macheps)\), then
\begin{equation}
\begin{split}
    \abs{\Delta(LU)^{(\mathtt{getrf})}_{rt}(n,n_b)}
    \leq{}&\bigO(\macheps)
    \norm{\hat L_{r,:}}\norm{\hat U_{:,t}}.
\end{split}
\end{equation}
\end{theorem}

\begin{proof}
For every preceding trailing update \(k=1,\ldots,s-1\), Theorem~\ref{lem:bound-XY-normwise},~\eqref{eq:lu-sum-bound}, and~\eqref{eq:lu-d-normwise} give
\begin{equation}
\begin{split}
    &\abs{\bigl(\Delta(L_{r,(k-1)n_b+1:kn_b}U_{(k-1)n_b+1:kn_b,t})^{(\mathtt{gemm})}(n_b)\bigr)^{(k)}} \\
    &\qquad\leq
    \bigl((1+\macheps)\varepsilon+\macheps\bigr)
    \norm{\hat L_{r,(k-1)n_b+1:kn_b}}
    \norm{\hat U_{(k-1)n_b+1:kn_b,t}}.
\end{split}
\end{equation}
Summing over \(k\) and applying Cauchy--Schwarz over the block index yields
\begin{equation*}
\begin{split}
    &\sum_{k=1}^{s-1}
    \abs{\bigl(\Delta(L_{r,(k-1)n_b+1:kn_b}U_{(k-1)n_b+1:kn_b,t})^{(\mathtt{gemm})}(n_b)\bigr)^{(k)}}\\
    &\qquad\leq
    \bigl((1+\macheps)\varepsilon+\macheps\bigr)
    \norm{\hat L_{r,:}}\norm{\hat U_{:,t}}.
\end{split}
\end{equation*}
Together with Corollary~\ref{cor:block-lu} and
\(
    \abs{\hat L_{rl}}\abs{\hat U_{lt}}
    \leq\norm{\hat L_{r,:}}\norm{\hat U_{:,t}},
\)
we obtain
\begin{equation*}
\begin{split}
    \abs{\Delta(LU)^{(\mathtt{getrf})}_{rt}(n,n_b)}
    \leq{}&\Bigl((1+\gamma_{s-1})\bigl((1+\macheps)\varepsilon+\macheps+\gamma_{n_b}\bigr)
    +\gamma_{s-1}\Bigr)
    \norm{\hat L_{r,:}}\norm{\hat U_{:,t}},
\end{split}
\end{equation*}
which proves~\eqref{eq:lu-gemm-normwise-final} by using \(1+\macheps\leq2\) and \(\macheps\leq\gamma_{n_b}\).
\end{proof}

Consequently, the normwise \(d\)-selection strategy developed for matrix multiplication can be used directly in block LU factorization.
This strategy requires essentially no additional work beyond the factorization itself, but it does not use the fact that the final LU backward error is measured relative to the complete row and column norms of the computed factors.

\subsubsection{An LU-aware adaptive selection}
\label{subsubsec:lu-normwise-adaptive}

We next develop a less restrictive normwise criterion by exploiting information already available from the computed parts of \(\hat L\) and \(\hat U\).
The error of an individual trailing update need not be small relative only to the row and column norms of the current block product, because the final LU backward error is measured relative to \(\norm{\hat L_{r,:}}\norm{\hat U_{:,s}}\).
At the \(i\)th trailing update, the quantities \(\norm{\hat L_{r,1:in_b}}\) and \(\norm{\hat U_{1:in_b,s}}\) are already available from the computed parts of the factors, and satisfy
\[
    \norm{\hat L_{r,(i-1)n_b+1:in_b}}\leq\norm{\hat L_{r,1:in_b}}\leq\norm{\hat L_{r,:}},
    \qquad
    \norm{\hat U_{(i-1)n_b+1:in_b,s}}\leq\norm{\hat U_{1:in_b,s}}\leq\norm{\hat U_{:,s}}.
\]

We therefore choose the smallest positive integer \(d_i\) such that
\begin{equation} \label{eq:lu-d-normwise-adaptive}
\begin{split}
    &4n_b(d_i+1)\theta^{d_i}
    \left(
    \max_r
    \frac{\norm{\hat L_{r,(i-1)n_b+1:in_b}}}
         {\norm{\hat L_{r,1:in_b}}}
    \right)
    \left(
    \max_s
    \frac{\norm{\hat U_{(i-1)n_b+1:in_b,s}}}
         {\norm{\hat U_{1:in_b,s}}}
    \right)
    \leq\frac{\varepsilon}{2},
\end{split}
\end{equation}
where a ratio \(0/0\) is defined to be zero.
Since both ratios in~\eqref{eq:lu-d-normwise-adaptive} are at most one, this criterion is never more restrictive than~\eqref{eq:lu-d-normwise}.
Let \(d_1\) be the smallest positive integer satisfying the fixed criterion~\eqref{eq:lu-d-normwise}.
Then \(d_1\) also satisfies~\eqref{eq:lu-d-normwise-adaptive} at every later iteration.
Therefore,
\[
    d_i\leq d_1,\qquad i=1,\ldots,p-1,
\]
and the first split count is the largest one required by the adaptive strategy.
At later updates, the accumulated prefix norms may make the criterion less restrictive and allow a smaller \(d_i\).

Indeed, by~\eqref{eq:bound-XY-trunc-component}, for every entry in the \(i\)th trailing update,
\begin{equation} \label{eq:lu-trunc-normwise-adaptive}
\begin{split}
    &\abs{\bigl(\Delta(
    L_{r,(i-1)n_b+1:in_b}
    U_{(i-1)n_b+1:in_b,s}
    )^{(\mathtt{trunc})}(d_i,n_b)\bigr)^{(i)}}\\
    &\qquad\leq
    4n_b(d_i+1)\theta^{d_i}
    \norm{\hat L_{r,(i-1)n_b+1:in_b}}
    \norm{\hat U_{(i-1)n_b+1:in_b,s}}\\
    &\qquad\leq
    \frac{\varepsilon}{2}
    \norm{\hat L_{r,1:in_b}}
    \norm{\hat U_{1:in_b,s}}\\
    &\qquad\leq
    \frac{\varepsilon}{2}
    \norm{\hat L_{r,:}}
    \norm{\hat U_{:,s}}.
\end{split}
\end{equation}
If \(\varepsilon\) satisfies~\eqref{eq:sum-half-eps} with \(d=d_1\), then the same condition holds for every \(d_i\leq d_1\).
Hence, the summation error is bounded by \(\varepsilon/2\) times the corresponding prefix-norm product.
Together with~\eqref{eq:lu-trunc-normwise-adaptive}, this gives
\begin{equation} \label{eq:lu-ozaki-normwise-adaptive-local}
\begin{split}
    &\abs{\bigl(\Delta(
    L_{r,(i-1)n_b+1:in_b}
    U_{(i-1)n_b+1:in_b,s}
    )^{(\mathtt{sum})}(d_i)\bigr)^{(i)}}\\
    &\quad+
    \abs{\bigl(\Delta(
    L_{r,(i-1)n_b+1:in_b}
    U_{(i-1)n_b+1:in_b,s}
    )^{(\mathtt{trunc})}(d_i,n_b)\bigr)^{(i)}}\\
    &\qquad\leq
    \varepsilon
    \norm{\hat L_{r,1:in_b}}
    \norm{\hat U_{1:in_b,s}}.
\end{split}
\end{equation}
The required prefix norms can be updated recursively from the current block row and column norms.
More precisely, after the current blocks have been computed, we use
\begin{equation} \label{eq:lu-prefix-norm-update}
\begin{split}
    \norm{\hat L_{r,1:in_b}}^2
    &=
    \norm{\hat L_{r,1:(i-1)n_b}}^2
    +
    \norm{\hat L_{r,(i-1)n_b+1:in_b}}^2,\\
    \norm{\hat U_{1:in_b,s}}^2
    &=
    \norm{\hat U_{1:(i-1)n_b,s}}^2
    +
    \norm{\hat U_{(i-1)n_b+1:in_b,s}}^2.
\end{split}
\end{equation}
Thus, only the squared norms of the newly computed block row and column segments need to be added at each iteration.
For the first trailing update, the prefix norms coincide with the current block row and column norms, so \(d_1\) is exactly the split count selected by the fixed criterion~\eqref{eq:lu-d-normwise}.
Because \(d_i\leq d_1\) for all later iterations, it is sufficient to check once whether \(d_1\leq d_{\max}\).
When this condition holds, every subsequent adaptive choice automatically lies in \(\{1,\ldots,d_1\}\), so no fallback to ordinary matrix multiplication is needed in the adaptive normwise algorithm.
The additional work therefore consists only of updating these vector norms, evaluating the ratios in~\eqref{eq:lu-d-normwise-adaptive}, and taking two maxima, and is negligible compared with the trailing matrix multiplication.
The resulting procedure is summarized in Algorithm~\ref{alg:block-LU-ozaki-normwise-adaptive}.

\begin{algorithm}[!tb]
\caption{Block LU with normwise adaptive fast Ozaki updates}
\label{alg:block-LU-ozaki-normwise-adaptive}
\begin{algorithmic}[1]
\REQUIRE A block size \(n_b\), a maximum number of splits \(d_{\max}\), a tolerance \(\varepsilon\), and \(\Acal\in\mathbb F^{m\times n}\), with \(m=qn_b\) and \(n=pn_b\).
\ENSURE A lower trapezoidal matrix \(\LL\) and an upper triangular matrix \(\UU\).
\STATE Choose \(d_1\) as the smallest positive integer satisfying~\eqref{eq:lu-d-normwise}.
\STATE Assume \(d_1\leq d_{\max}\).
\FOR{\(i=1:p-1\)}
    \STATE Compute the LU factorization \(\Acal_{i:q,i}=\LL_{i:q,i}\UU_{i,i}\).
    \STATE Compute \(\UU_{i,i+1:p}\) from \(\LL_{i,i}\UU_{i,i+1:p}=\Acal_{i,i+1:p}\).
    \STATE Update the prefix row and column norms using~\eqref{eq:lu-prefix-norm-update}.
    \IF{\(i>1\)}
        \STATE Choose the smallest \(d_i\in\{1,\ldots,d_1\}\) satisfying~\eqref{eq:lu-d-normwise-adaptive}.
    \ENDIF
    \STATE Compute \(\LL_{i+1:q,i}\UU_{i,i+1:p}\) by the fast Ozaki scheme with \(d_i\) splits and update the trailing submatrix.
\ENDFOR
\STATE Compute the LU factorization \(\Acal_{p:q,p}=\LL_{p:q,p}\UU_{p,p}\).
\end{algorithmic}
\end{algorithm}

We then show the normwise stability of Algorithm~\ref{alg:block-LU-ozaki-normwise-adaptive} in Theorem~\ref{thm:lu-normwise-adaptive-d}.

\begin{theorem} \label{thm:lu-normwise-adaptive-d}
Consider Algorithm~\ref{alg:block-LU-ozaki-normwise-adaptive}.
Assume that \(d_1\leq d_{\max}\) and that \(\varepsilon\) satisfies~\eqref{eq:sum-half-eps} with \(d=d_1\).
Then the adaptive selection is well defined, with \(d_i\leq d_1\) for every \(i=1,\ldots,p-1\).
Moreover, there exists \(\Delta(LU)^{(\mathtt{getrf})}(n,n_b)=\Delta(\LL\UU)^{(\mathtt{getrf})}(n,n_b)\) such that~\eqref{eq:blocklu:getrf} holds with, for any scalar entry \((r,t)\) in block \((i,j)\), let \(s=\min\{i,j\}\),
\begin{equation} \label{eq:lu-gemm-normwise-adaptive-final}
\begin{split}
    \abs{\Delta(LU)^{(\mathtt{getrf})}_{rt}(n,n_b)}
    \leq{}&
    2(1+\gamma_{s-1})
    \bigl((s-1)\varepsilon+\gamma_{n_b}+\gamma_{s-1}\bigr)
    \norm{\hat L_{r,:}}\norm{\hat U_{:,t}}.
\end{split}
\end{equation}
In particular, if \(\varepsilon=\bigO(\macheps)\), then
\[
    \abs{\Delta(LU)^{(\mathtt{getrf})}_{rt}(n,n_b)}
    \leq
    \bigO(\macheps)\norm{\hat L_{r,:}}\norm{\hat U_{:,t}}.
\]
\end{theorem}

\begin{proof}
By construction, \(d_1\) satisfies the fixed criterion~\eqref{eq:lu-d-normwise}.
Since both ratios in~\eqref{eq:lu-d-normwise-adaptive} are at most one, the same \(d_1\) satisfies the adaptive criterion at every later iteration.
Hence, the set \(\{1,\ldots,d_1\}\) always contains an admissible split count, so the adaptive choice is well defined and \(d_i\leq d_1\).

Consider a preceding trailing update \(k=1,\ldots,s-1\).
Since \(\varepsilon\) satisfies~\eqref{eq:sum-half-eps} with \(d=d_1\) and \(d_k\leq d_1\),~\eqref{eq:lu-ozaki-normwise-adaptive-local} and~\eqref{eq:bound-DLU} give
\begin{equation} \label{eq:lu-gemm-entry-normwise-adaptive}
\begin{split}
    &\abs{\bigl(\Delta(
    L_{r,(k-1)n_b+1:kn_b}
    U_{(k-1)n_b+1:kn_b,t}
    )^{(\mathtt{gemm})}(n_b)\bigr)^{(k)}}\\
    &\qquad\leq
    (1+\macheps)\varepsilon
    \norm{\hat L_{r,1:kn_b}}
    \norm{\hat U_{1:kn_b,t}}
    +
    \macheps
    \norm{\hat L_{r,(k-1)n_b+1:kn_b}}
    \norm{\hat U_{(k-1)n_b+1:kn_b,s}}.
\end{split}
\end{equation}
The first norm product is bounded by
\(\norm{\hat L_{r,:}}\norm{\hat U_{:,t}}\).
Therefore, summing the first terms in~\eqref{eq:lu-gemm-entry-normwise-adaptive} over \(k=1,\ldots,s-1\) gives \((s-1)(1+\macheps)\varepsilon\norm{\hat L_{r,:}}\norm{\hat U_{:,t}}\).
For the remaining term, applying the Cauchy--Schwarz inequality over the block index gives
\begin{equation*}
\begin{split}
    &\sum_{k=1}^{s-1}
    \norm{\hat L_{r,(k-1)n_b+1:kn_b}}
    \norm{\hat U_{(k-1)n_b+1:kn_b,t}}\\
    &\qquad\leq
    \norm{\hat L_{r,1:(s-1)n_b}}
    \norm{\hat U_{1:(s-1)n_b,t}}
    \leq
    \norm{\hat L_{r,:}}\norm{\hat U_{:,t}}.
\end{split}
\end{equation*}
Hence,
\begin{equation} \label{eq:lu-gemm-sum-normwise-adaptive}
\begin{split}
    &\sum_{k=1}^{s-1}
    \abs{\bigl(\Delta(
    L_{r,(k-1)n_b+1:kn_b}
    U_{(k-1)n_b+1:kn_b,t}
    )^{(\mathtt{gemm})}(n_b)\bigr)^{(k)}}\\
    &\qquad\leq
    \bigl((s-1)(1+\macheps)\varepsilon+\macheps\bigr)
    \norm{\hat L_{r,:}}\norm{\hat U_{:,t}}.
\end{split}
\end{equation}
Applying Corollary~\ref{cor:block-lu} entrywise and using
\(
    \abs{\hat L_{rl}}\abs{\hat U_{lt}}
    \leq\norm{\hat L_{r,:}}\norm{\hat U_{:,t}},
\)
we obtain
\begin{equation*}
\begin{split}
    \abs{\Delta(LU)^{(\mathtt{getrf})}_{rt}(n,n_b)}
    \leq{}&\Bigl((1+\gamma_{s-1})\bigl((s-1)(1+\macheps)\varepsilon+\macheps+\gamma_{n_b}\bigr)
    +\gamma_{s-1}\Bigr)
    \norm{\hat L_{r,:}}\norm{\hat U_{:,t}},
\end{split}
\end{equation*}
which proves~\eqref{eq:lu-gemm-normwise-adaptive-final} by using \(1+\macheps\leq2\) and \(\macheps\leq\gamma_{n_b}\).
\end{proof}

\subsection{Componentwise stability}
\label{subsec:lu-componentwise}
Throughout this subsection, the indices \(i\) and \(k\) are block indices.
The basic idea is analogous to the LU-aware normwise strategy in Section~\ref{subsubsec:lu-normwise-adaptive}: the accuracy required of an individual trailing matrix product can be relaxed by using information accumulated in the factorization.
For componentwise stability, however, directly applying the strategy in Section~\ref{subsec:alo-component} to every trailing update can be unnecessarily restrictive, because the lower-reference factor used for an individual matrix product can be very small or even zero.
Corollary~\ref{cor:block-lu} shows that the errors from different trailing updates enter the final LU backward error additively.
We therefore construct a componentwise reference matrix \(\GG\) from the matrix products already generated during the factorization.

Let \(\GG\) be a nonnegative matrix initialized to zero.
At the first iteration, no previous trailing update is available and hence \(\GG=0\).
We therefore select \(d\) for the first trailing product by the componentwise matrix-multiplication strategy in Section~\ref{subsec:alo-component}.
If the selected \(d\leq d_{\max}\), the fast Ozaki scheme is used.
Otherwise, the first trailing product is evaluated by ordinary matrix multiplication.
For an accepted Ozaki update, the componentwise selection rule in Section~\ref{subsec:alo-component}, together with~\eqref{eq:sum-half-eps}, gives
\begin{equation} \label{eq:lu-first-componentwise-generic}
\begin{split}
    &\abs{\bigl(\Delta(\LL_{2:q,1}\UU_{1,2:p})^{(\mathtt{sum})}(d)\bigr)^{(1)}}
    +\abs{\bigl(\Delta(\LL_{2:q,1}\UU_{1,2:p})^{(\mathtt{trunc})}(d,n_b)\bigr)^{(1)}}
    \leq\varepsilon\abs{\LLhat_{2:q,1}}\abs{\UUhat_{1,2:p}}.
\end{split}
\end{equation}

After the first trailing product has been evaluated, it provides the first contribution to \(\GG\).
If the fast Ozaki scheme is used, we set
\[
    \GG_{2:q,2:p}=\abs{\LLhat_{2:q,1}^{(1)}}\abs{\UUhat_{1,2:p}^{(1)}}.
\]
If ordinary matrix multiplication is used, denote the computed product by \(\CC^{(1)}\) and set
\[
    \GG_{2:q,2:p}=\abs{\CC^{(1)}}.
\]
For every later iteration \(i>1\), \(d\) is selected from the current \(\GG\) by requiring
\begin{equation} \label{eq:lu-d-componentwise}
    n_b(d+1)\theta^d\widetilde{x}^{(i)}\bigl(\widetilde{y}^{(i)}\bigr)\trans\leq\varepsilon\GG_{i+1:q,i+1:p}
\end{equation}
componentwise, where
\begin{equation} \label{eq:lu-scaling-factors}
    \widetilde{x}^{(i)}_r=2^{\left\lceil\log_2\max_{l}\abs{(\LLhat_{i+1:q,i})_{rl}}\right\rceil},\qquad
    \widetilde{y}^{(i)}_s=2^{\left\lceil\log_2\max_{l}\abs{(\UUhat_{i,i+1:p})_{ls}}\right\rceil}.
\end{equation}
If the selected \(d\leq d_{\max}\), the fast Ozaki scheme is used.
Otherwise, the trailing product is evaluated by ordinary matrix multiplication.

For every accepted Ozaki update, the corresponding first-split product \(\abs{\LLhat_{i+1:q,i}^{(1)}}\abs{\UUhat_{i,i+1:p}^{(1)}}\) is added to \(\GG\).
The first splits are already represented in the low precision.
Taking their absolute values introduces no rounding error, and the choice of the splitting parameter in the fast Ozaki scheme guarantees that this product is computed exactly on the Tensor Core.
By~\eqref{eq:sum-splits-X-Y} in Lemma~\ref{lem:bound-sum-X},
\[
    \abs{\LLhat_{i+1:q,i}^{(1)}}\leq\sum_{r\geq1}\abs{\LLhat_{i+1:q,i}^{(r)}}\leq\kappa_{\rm sp}\abs{\LLhat_{i+1:q,i}},
    \qquad
    \abs{\UUhat_{i,i+1:p}^{(1)}}\leq\sum_{s\geq1}\abs{\UUhat_{i,i+1:p}^{(s)}}\leq\kappa_{\rm sp}\abs{\UUhat_{i,i+1:p}},
\]
and hence
\begin{equation} \label{eq:lu-G-first-product-bound}
    \abs{\LLhat_{i+1:q,i}^{(1)}}\abs{\UUhat_{i,i+1:p}^{(1)}}
    \leq\kappa_{\rm sp}^2\abs{\LLhat_{i+1:q,i}}\abs{\UUhat_{i,i+1:p}}.
\end{equation}

For an update evaluated by ordinary matrix multiplication, denote the computed product by \(\CC^{(i)}\).
Its absolute value is used to update \(\GG\).
The standard componentwise matrix-multiplication bound gives
\begin{equation} \label{eq:lu-G-ordinary-product-bound}
    \abs{\CC^{(i)}}
    \leq(1+\gamma_{n_b})
    \abs{\LLhat_{i+1:q,i}}\abs{\UUhat_{i,i+1:p}}.
\end{equation}
Thus, after the first iteration, every trailing update provides a usable LU-scale reference for the subsequent \(d\)-selection.

The resulting procedure is summarized in Algorithm~\ref{alg:block-LU-ozaki-componentwise}.

\begin{algorithm}[!tb]
\caption{Block LU with componentwise adaptive fast Ozaki updates}
\label{alg:block-LU-ozaki-componentwise}
\begin{algorithmic}[1]
\REQUIRE A block size \(n_b\), a maximum number of splits \(d_{\max}\), a tolerance \(\varepsilon\), and \(\Acal\in\mathbb F^{m\times n}\), with \(m=qn_b\) and \(n=pn_b\).
\ENSURE A lower trapezoidal matrix \(\LL\) and an upper triangular matrix \(\UU\).
\STATE Set \(\GG=0\).
\FOR{\(i=1:p-1\)}
    \STATE Compute the LU factorization \(\Acal_{i:q,i}=\LL_{i:q,i}\UU_{i,i}\).
    \STATE Compute \(\UU_{i,i+1:p}\) from \(\LL_{i,i}\UU_{i,i+1:p}=\Acal_{i,i+1:p}\).
    \IF{\(i=1\)}
        \STATE Select \(d\) using the componentwise matrix-multiplication strategy in Section~\ref{subsec:alo-component}.
    \ELSE
        \STATE Find the smallest positive integer \(d\) satisfying~\eqref{eq:lu-d-componentwise}; set \(d=\infty\) if no finite \(d\) exists.
    \ENDIF
    \IF{\(d\leq d_{\max}\)}
        \STATE Compute \(\LL_{i+1:q,i}\UU_{i,i+1:p}\) by the fast Ozaki scheme with \(d\) splits and update the trailing submatrix.
        \STATE Compute \(\abs{\LL_{i+1:q,i}^{(1)}}\abs{\UU_{i,i+1:p}^{(1)}}\) exactly on the Tensor Core.
        \STATE Update \(\GG_{i+1:q,i+1:p}=\GG_{i+1:q,i+1:p}+\abs{\LL_{i+1:q,i}^{(1)}}\abs{\UU_{i,i+1:p}^{(1)}}\).
    \ELSE
        \STATE Compute \(\CC^{(i)}=\LL_{i+1:q,i}\UU_{i,i+1:p}\) by ordinary matrix multiplication and update the trailing submatrix.
        \STATE Update \(\GG_{i+1:q,i+1:p}=\GG_{i+1:q,i+1:p}+\abs{\CC^{(i)}}\).
    \ENDIF
\ENDFOR
\STATE Compute the LU factorization \(\Acal_{p:q,p}=\LL_{p:q,p}\UU_{p,p}\).
\end{algorithmic}
\end{algorithm}

\begin{lemma} \label{lem:lu-G-rounding}
Assume that no overflow or underflow occurs.
Suppose that \(t\geq1\) trailing updates have contributed to the matrix \(\GG\) used by Algorithm~\ref{alg:block-LU-ozaki-componentwise}.
Then the computed matrix \(\GGhat\) satisfies
\begin{equation} \label{eq:lu-G-rounding}
    0\leq\GGhat
    \leq
    (1+\gamma_{t-1})
    \max\{\kappa_{\rm sp}^2,1+\gamma_{n_b}\}
    \abs{\LLhat}\abs{\UUhat}.
\end{equation}
In particular, since \(t\leq p-1\),
\begin{equation} \label{eq:lu-G-rounding-p}
    \GGhat
    \leq
    (1+\gamma_{p-2})
    \max\{\kappa_{\rm sp}^2,1+\gamma_{n_b}\}
    \abs{\LLhat}\abs{\UUhat}.
\end{equation}
\end{lemma}

\begin{proof}
An accepted Ozaki update contributes its first-split product, which is bounded by~\eqref{eq:lu-G-first-product-bound}.
An update evaluated by ordinary matrix multiplication contributes \(\abs{\CC^{(i)}}\), which is bounded by~\eqref{eq:lu-G-ordinary-product-bound}.
Hence, the sum of the exact nonnegative contributions embedded in their corresponding block positions is bounded componentwise by
\[
    \max\{\kappa_{\rm sp}^2,1+\gamma_{n_b}\}
    \sum_k\abs{\LLhat_{:,k}}\abs{\UUhat_{k,:}}
    \leq\max\{\kappa_{\rm sp}^2,1+\gamma_{n_b}\}\abs{\LLhat}\abs{\UUhat}.
\]
Taking absolute values introduces no rounding error.
The first contribution is assigned directly to \(\GG\), and the remaining \(t-1\) contributions are accumulated recursively in precision \(\macheps\).
The standard rounding error bound for recursive summation therefore gives~\eqref{eq:lu-G-rounding}.
\end{proof}

\begin{theorem} \label{thm:lu-componentwise-d}
Consider Algorithm~\ref{alg:block-LU-ozaki-componentwise}.
Assume that \(\varepsilon\) satisfies~\eqref{eq:sum-half-eps} with \(d=d_{\max}\).
At the first iteration, assume that \(d\) is selected by the componentwise matrix-multiplication strategy in Section~\ref{subsec:alo-component}.
At every later accepted fast Ozaki update, assume that \(d\) satisfies~\eqref{eq:lu-d-componentwise}.
Also assume that every ordinary matrix product \(\CC=XY\) satisfies
\[
    \abs{\CC-XY}
    \leq
    \gamma_{n_b}\abs{X}\abs{Y}.
\]
Then, for the \(i\)th block column,
\begin{equation} \label{eq:lu-componentwise-final}
\begin{split}
    \abs{\Delta(\LL\UU)^{(\mathtt{getrf})}_{i:q,i}(n,n_b)}
    \leq{}&2(1+\gamma_{i-1})\Bigl(\varepsilon
    +\gamma_{d_{\max}(d_{\max}+1)/2}\kappa_{\rm sp}^2
    +2\gamma_{n_b}+\gamma_{i-1}\\
    &+\max\{i-2,0\}\varepsilon(1+\gamma_{p-2})
    \max\{\kappa_{\rm sp}^2,1+\gamma_{n_b}\}\Bigr)
    \abs{\LLhat_{i:q,1:i}}\abs{\UUhat_{1:i,i}}.
\end{split}
\end{equation}
In particular, if \(\varepsilon=\bigO(\macheps)\) and \(n_b\leq 2^{21}\), then \(\kappa_{\rm sp}\leq 7\) and
\begin{equation}
\begin{split}
    \abs{\Delta(LU)^{(\mathtt{getrf})}(n,n_b)}
    \leq{}&\bigO(\macheps)\abs{\LLhat}\abs{\UUhat}.
\end{split}
\end{equation}
\end{theorem}

\begin{proof}
By Lemma~\ref{lem:lu-G-rounding}, after the first trailing update,
\begin{equation} \label{eq:lu-G-bound-theorem}
    \GGhat
    \leq
    (1+\gamma_{p-2})
    \max\{\kappa_{\rm sp}^2,1+\gamma_{n_b}\}
    \abs{\LLhat}\abs{\UUhat}.
\end{equation}
For the first trailing update, if the fast Ozaki scheme is used,~\eqref{eq:lu-first-componentwise-generic} and~\eqref{eq:bound-DLU} give
\[
    \abs{\bigl(\Delta(\LL_{2:q,1}\UU_{1,2:p})^{(\mathtt{gemm})}(n_b)\bigr)^{(1)}}
    \leq
    \bigl((1+\macheps)\varepsilon+\macheps\bigr)
    \abs{\LLhat_{2:q,1}}\abs{\UUhat_{1,2:p}}.
\]
If the first update is evaluated by ordinary matrix multiplication, the standard product bound and~\eqref{eq:bound-DLU} give
\[
    \abs{\bigl(\Delta(\LL_{2:q,1}\UU_{1,2:p})^{(\mathtt{gemm})}(n_b)\bigr)^{(1)}}
    \leq
    2\gamma_{n_b}
    \abs{\LLhat_{2:q,1}}\abs{\UUhat_{1,2:p}}.
\]
Hence, in either case,
\begin{equation} \label{eq:lu-first-gemm-componentwise}
    \abs{\bigl(\Delta(\LL_{2:q,1}\UU_{1,2:p})^{(\mathtt{gemm})}(n_b)\bigr)^{(1)}}
    \leq
    2(\varepsilon+\gamma_{n_b})
    \abs{\LLhat_{2:q,1}}\abs{\UUhat_{1,2:p}}.
\end{equation}

For every accepted Ozaki update with \(k>1\),~\eqref{eq:bound-XY-trunc-component} and~\eqref{eq:lu-d-componentwise} give
\begin{equation} \label{eq:lu-trunc-componentwise-G}
\begin{split}
    \abs{\bigl(\Delta(\LL_{k+1:q,k}\UU_{k,k+1:p})^{(\mathtt{trunc})}(d,n_b)\bigr)^{(k)}}
    &\leq\varepsilon\GGhat_{k+1:q,k+1:p}\\
    &\leq\varepsilon(1+\gamma_{p-2})
    \max\{\kappa_{\rm sp}^2,1+\gamma_{n_b}\}
    \bigl(\abs{\LLhat}\abs{\UUhat}\bigr)_{k+1:q,k+1:p}.
\end{split}
\end{equation}
Moreover, since every accepted \(d\leq d_{\max}\),
\begin{equation} \label{eq:lu-sum-componentwise-dmax}
    \abs{\bigl(\Delta(\LL_{k+1:q,k}\UU_{k,k+1:p})^{(\mathtt{sum})}(d)\bigr)^{(k)}}
    \leq
    \gamma_{d_{\max}(d_{\max}+1)/2}\kappa_{\rm sp}^2
    \abs{\LLhat_{k+1:q,k}}\abs{\UUhat_{k,k+1:p}}.
\end{equation}
For an ordinary matrix multiplication,
\begin{equation} \label{eq:lu-gemm-componentwise-ordinary}
\begin{split}
    \abs{\bigl(\Delta(\LL_{k+1:q,k}\UU_{k,k+1:p})^{(\mathtt{gemm})}(n_b)\bigr)^{(k)}}
    &\leq
    \bigl((1+\macheps)\gamma_{n_b}+\macheps\bigr)
    \abs{\LLhat_{k+1:q,k}}\abs{\UUhat_{k,k+1:p}}\\
    &\leq
    2\gamma_{n_b}
    \abs{\LLhat_{k+1:q,k}}\abs{\UUhat_{k,k+1:p}}.
\end{split}
\end{equation}

Combining~\eqref{eq:lu-first-gemm-componentwise} with the bounds for the later updates, and using the fact that the local block products form disjoint terms in
\(\abs{\LLhat_{i:q,1:i-1}}\abs{\UUhat_{1:i-1,i}}\), gives
\begin{equation} \label{eq:lu-gemm-sum-componentwise}
\begin{split}
    &\sum_{k=1}^{i-1}
    \abs{\bigl(\Delta(\LL_{i:q,k}\UU_{k,i})^{(\mathtt{gemm})}(n_b)\bigr)^{(k)}}\\
    &\quad\leq
    2\Bigl(
    \varepsilon
    +\gamma_{d_{\max}(d_{\max}+1)/2}\kappa_{\rm sp}^2
    +\gamma_{n_b}
    +\max\{i-2,0\}\varepsilon(1+\gamma_{p-2})
    \max\{\kappa_{\rm sp}^2,1+\gamma_{n_b}\}
    \Bigr)\\
    &\qquad\cdot
    \abs{\LLhat_{i:q,1:i}}\abs{\UUhat_{1:i,i}}.
\end{split}
\end{equation}
Applying Corollary~\ref{cor:block-lu} to~\eqref{eq:lu-gemm-sum-componentwise} and using \(\abs{\LLhat_{i:q,1:i-1}}\abs{\UUhat_{1:i-1,i}}\leq \abs{\LLhat_{i:q,1:i}}\abs{\UUhat_{1:i,i}}\) yields~\eqref{eq:lu-componentwise-final}.
\end{proof}

\subsection{Computational costs of \(d\)-selection strategies}
\label{subsec:lu-additional-cost}

We briefly compare the additional costs of the three strategies for selecting \(d\).
Let \(c_d\) denote the cost of one scalar search for \(d\), \(c_{\rm norm}\) the cost constant for updating the row and column norms, \(c_{\rm elem}\) the cost constant for elementwise operations on the trailing matrices, \(c_{\rm sum}\) the cost constant for rescaling and accumulating one split product in precision \(\macheps\), \(c_{\rm lu}\) the leading working precision cost constant of the block LU factorization, and \(c_{\rm gemm}\) the cost constant of a Tensor Core matrix multiplication with low precision inputs and medium precision accumulation.

\begin{table}[htbp]
\centering
\caption{Additional arithmetic costs of the three \(d\)-selection strategies and the leading cost of block LU with Ozaki matrix multiplication.
Here we assume that the trailing updates of block LU are all performed by the fast Ozaki Scheme.}
\label{tab:lu-d-selection-cost}
\begin{tabular}{l|cc}
\hline
&
\(\macheps\)
&
\(\machepsmedm/\machepslow\) on Tensor Cores\\
\hline
\begin{tabular}[c]{@{}l@{}}Fixed norm\\-wise selection\end{tabular}
& \(c_d\)
& \(0\)\\
\begin{tabular}[c]{@{}l@{}}Adaptive norm\\-wise selection\end{tabular}
& \(c_d\dfrac{n}{n_b}+c_{\rm norm}mn\)
& \(0\)\\
\begin{tabular}[c]{@{}l@{}l@{}}Adaptive\\ componentwise\\ selection\end{tabular}
& \(c_d\dfrac{n}{n_b}+c_{\rm elem}\dfrac{mn^2}{n_b}\)
& \(c_{\rm gemm}mn^2\)\\
\hline
\begin{tabular}[c]{@{}l@{}l@{}}Block LU + \\Ozaki matrix \\ multiplication\end{tabular}
& \(c_{\rm lu}mn^2+
c_{\rm sum}\left(\displaystyle\frac{1}{p}\sum_{i=1}^{p}\frac{d_i(d_i+1)}{2}\right)\dfrac{mn^2}{n_b}\)
& \(c_{\rm gemm}
\left(\displaystyle\frac{1}{p}\sum_{i=1}^{p}\frac{d_i(d_i+1)}{2}\right)mn^2\)\\
\hline
\end{tabular}
\end{table}

The fixed normwise selection has the smallest additional cost.
The value of \(d\) is selected only once before the factorization, so the cost of searching for \(d\) is only \(c_d\), which is usually \(\bigO(1)\), and the same value is then used at every iteration.
By contrast, both adaptive strategies select \(d_i\) at every block iteration, so the scalar-selection work is bounded by \(c_dn/n_b\).
The adaptive normwise selection also requires extra work to update the prefix row and column norms, but it may select smaller values of \(d_i\) in later iterations.

The additional costs of all three selection strategies are small compared with the overall cost of block LU with Ozaki matrix multiplication.
In particular,
\[
    c_{\rm elem}\ll c_{\rm lu}n_b+c_{\rm sum}
    \left(\frac{1}{p}\sum_{i=1}^{p}\frac{d_i(d_i+1)}{2}\right),
    \qquad
    1\ll\frac{1}{p}\sum_{i=1}^{p}\frac{d_i(d_i+1)}{2}.
\]
Hence, the extra work required for selecting \(d\) is expected to be modest relative to the complete factorization.

\section{Numerical experiments}
\label{sec:experiments}

We present numerical experiments for the normwise and componentwise choices of the number of splits in the Ozaki matrix multiplication.
Throughout the experiments, $\macheps$, $\machepsmedm$, and $\machepslow$ denote the unit roundoffs of double, single, and half precision, respectively.
Since our purpose here is only to examine the accuracy of the proposed strategies to select $d$, all experiments are carried out in MATLAB rather than on actual Tensor Cores.
In our MATLAB implementation, we select the splitting parameter to satisfy the Ozaki scheme’s exact-accumulation condition.
The split matrix products are computed using single-precision matrix multiplication, which is exact under the splitting assumptions, and the resulting terms are accumulated in double precision.
The normwise and componentwise backward errors are measured, respectively, by
\[
\max_{r,t}\frac{\abs{(P A-\widehat L\widehat U)_{rt}}}{\norm{\widehat L_{r,:}}\norm{\widehat U_{:,t}}},
\qquad
\max_{r,t}\frac{\abs{(P A-\widehat L\widehat U)_{rt}}}{(\abs{\widehat L}\abs{\widehat U})_{rt}},
\]
where \(P\) denote the permutation matrix generated by partial pivoting.

\subsection{A constructed example}
\label{subsec:experiment-constructed}

We first construct a small example to illustrate the difference between the normwise and componentwise choices of $d$.
We take $m=n=8$ and $n_b=4$, so the factorization contains only one trailing block update.
Consequently, the normwise adaptive criterion reduces to its first-update criterion, while the componentwise strategy uses the componentwise matrix-multiplication criterion of Section~\ref{subsec:alo-component}.
We construct
\[
\LL=\begin{bmatrix}
I&0\\
\LL_{21}&I
\end{bmatrix},
\qquad
\UU=\begin{bmatrix}
I&\UU_{12}\\
0&I
\end{bmatrix},
\qquad
A=\LL\UU,
\]
where $I$ denotes the $4\times4$ identity matrix.
Let \(\delta=2^{-26}\) and \(\tau=2^{-12}\).
The two off-diagonal blocks are chosen as
\[
\LL_{21}=
\begin{bmatrix}
1+\tau&\delta(1+\tau)&\delta^2&\delta^3\\
\delta^3&1+3\tau&\delta(1+\tau)&\delta^2\\
\delta^2&\delta^3&1+5\tau&\delta(1+\tau)\\
\delta(1+\tau)&\delta^2&\delta^3&1+7\tau
\end{bmatrix}
\]
and
\[
\UU_{12}=
\begin{bmatrix}
\delta(1+2\tau)&\delta^2&\delta^3&1+2\tau\\
1+4\tau&\delta(1+4\tau)&\delta^2&\delta^3\\
\delta^3&1+6\tau&\delta(1+6\tau)&\delta^2\\
\delta^2&\delta^3&1+8\tau&\delta(1+8\tau)
\end{bmatrix}.
\]
For this example, $\kappa_2(A)=6.8676$, so the difference observed below is not caused by ill-conditioning of $A$.
The construction creates entries for which the normwise scale
$\norm{(\LL_{21})_{r,:}}\norm{(\UU_{12})_{:,t}}$
is much larger than the corresponding componentwise scale
$(\abs{\LL_{21}}\abs{\UU_{12}})_{rt}$.
For example, the largest ratio between these two scales is approximately $7.57\times10^{22}$.

Table~\ref{tab:constructed-example} reports the results with $\varepsilon=\macheps$ and $d_{\max}=30$.
The normwise criterion selects $d=6$.
Although the resulting normwise backward error is very small, the componentwise backward error is approximately $3.3\times10^{-1}$.
In contrast, the componentwise criterion selects $d=20$ and gives zero backward error in this experiment.
This example therefore demonstrates that a normwise choice of $d$ does not in general guarantee a small componentwise backward error.

\begin{table}[htbp]
\centering
\caption{Results for the constructed example with $m=n=8$, $n_b=4$ to compare normwise and componentwise \(d\)-selection strategies.
Here we use double precision to measure the backward errors.}
\label{tab:constructed-example}
\begin{tabular}{l|ccc}
\hline
$d$-selection & $d$ & Normwise backward error & Componentwise backward error\\
\hline
Normwise      & $6$  & $2.6447\times10^{-23}$ & $3.3328\times10^{-1}$\\
Componentwise & $20$ & $0$                     & $0$\\
\hline
\end{tabular}
\end{table}

\subsection{Tests with random matrices}
\label{subsec:experiment-gallery}

We next test the two adaptive choices of \(d\) on a broader collection of matrices generated by MATLAB's \texttt{gallery(`randsvd')} function.
We take \(m=n=128\), \(n_b=32\), and \(\varepsilon=n\macheps\).
The prescribed condition numbers are \(10^2\), \(10^4\), \(10^6\), \(10^8\), \(10^{10}\), and \(10^{12}\).
For each condition number, we use the five singular value distributions corresponding to \texttt{mode}\(=1,\ldots,5\), giving a total of 30 test matrices.
Since \(m=n=128\) and \(n_b=32\), there are three trailing block updates in each factorization.
For the componentwise strategy we set \(d_{\max}=20\).

\begin{figure}[htbp]
\centering
\includegraphics[width=\textwidth]{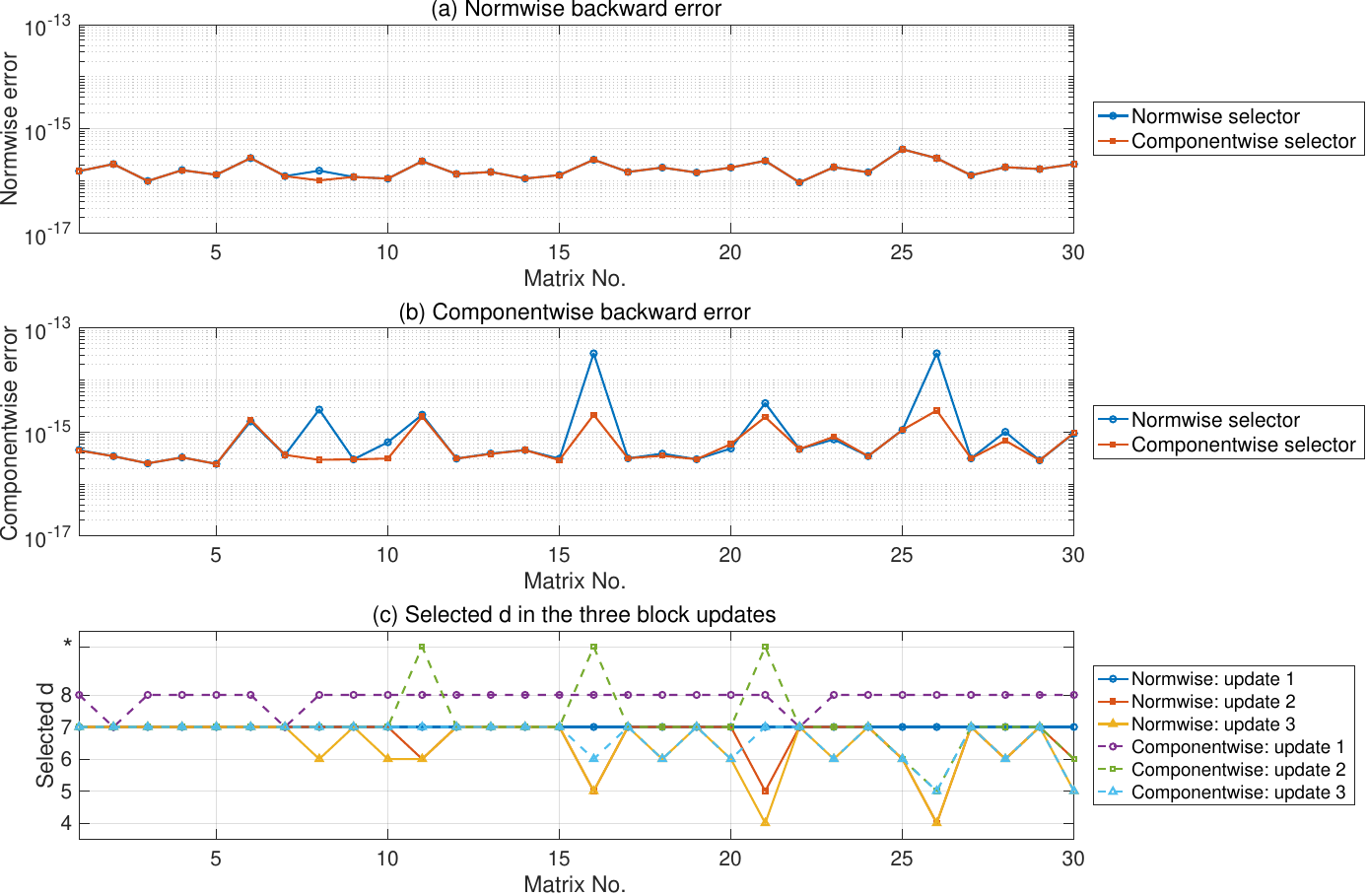}
\caption{Results for the 30 matrices generated by \texttt{gallery('randsvd')} with \(m=n=128\), \(n_b=32\), and \(\varepsilon=n\macheps\).
The first two panels show the final normwise and componentwise backward errors, respectively.
The third panel shows the values of \(d\) selected in the three trailing block updates for both adaptive strategies.
The symbol ``\(\star\)'' indicates an ordinary matrix multiplication because no \(d\leq d_{\max}=20\) satisfies the componentwise criterion.}
\label{fig:gallery-experiments}
\end{figure}

Figure~\ref{fig:gallery-experiments} shows the final normwise and componentwise backward errors and the selected values of \(d\) in the three trailing block updates.
The two strategies give almost indistinguishable normwise backward errors for all 30 matrices, with the \(\bigO(10^{-16})\) level.
Their componentwise behavior is more clearly different.
The normwise strategy usually also gives a small componentwise backward error, but several matrices show noticeably larger errors than those obtained with the componentwise strategy.
The difference is especially visible for a few test matrices for which the normwise componentwise error rises to around \(10^{-14}\), whereas the componentwise strategy remains around \(10^{-15}\).
Thus, the constructed example in Section~\ref{subsec:experiment-constructed} is deliberately extreme, but the advantage of using componentwise information is still visible for some of the random test matrices.

The third panel also illustrates the adaptivity of both selection strategies.
For these tests, the first normwise update uses \(d_1=7\), while later normwise updates can use smaller values, in some cases as small as \(d=4\).
This confirms that the LU-aware normwise criterion can reduce the number of split matrix products compared with using the first value \(d_1\) throughout the factorization.
Since an Ozaki update with \(d\) splits requires \(d(d+1)/2\) low precision matrix multiplications, reducing \(d\) from \(7\) to \(5\) decreases this number from \(28\) to \(15\) (around \(46\%\) reduction).
The componentwise strategy generally selects \(d=8\) for the first update and then adapts the later values according to the accumulated matrix \(\GG\).
The later componentwise choices are often between \(5\) and \(7\), although a few updates do not satisfy the criterion for any \(d\leq d_{\max}\) and therefore fall back to ordinary matrix multiplication, indicated by ``\(\star\)'' in the figure.

Overall, the experiments show that the normwise adaptive strategy can reduce the split depth at essentially negligible selection cost, while the componentwise strategy uses the evolving componentwise scale of the factorization to obtain better componentwise accuracy when needed.
The latter may require a larger split depth, or occasionally an ordinary matrix multiplication, but it provides a more reliable control of the componentwise backward error.

\section{Conclusion}
\label{sec:conclusion}

In this paper, we have studied how to select the number of splits \(d\) in the Ozaki matrix multiplication according to the accuracy required by the computation.
Since \(d\) controls both the number of split matrix products and the truncation error, its choice directly determines the tradeoff between computational work and accuracy.

For the Ozaki matrix multiplication itself, we derived normwise and componentwise criteria for selecting \(d\).
These criteria are obtained by controlling the truncation error together with the rounding error from accumulating the retained split products.
They therefore provide a direct way to determine how many split products are needed for a prescribed accuracy scale.

We then incorporated the same idea into block LU factorization, where the Ozaki scheme is used for the trailing matrix updates.
For normwise stability, we considered both a fixed strategy and an LU-aware adaptive strategy.
The fixed strategy selects \(d\) only once and has essentially no additional selection cost.
The adaptive strategy uses prefix row and column norms of the computed factors and can reduce the split depth in later updates.
In particular, the first split depth \(d_1\) is the largest one required by the adaptive normwise strategy, so all subsequent selections satisfy \(d_i \leq d_1\).
For componentwise stability, the first trailing update is selected by the componentwise matrix-multiplication criterion.
This first update initializes a nonnegative reference matrix \(\GG\), which is then accumulated throughout the factorization and used to select \(d\) in the later updates.
If no \(d\leq d_{\max}\) satisfies the componentwise criterion, the corresponding update is evaluated by ordinary matrix multiplication.
The resulting strategy therefore adapts both the split depth and, when necessary, the matrix-multiplication kernel to the componentwise scale of the factorization.

The choice between the normwise and componentwise strategies should therefore be guided by the accuracy required by the final application.
If normwise stability is sufficient, the normwise strategy is preferable: it requires less additional work for selecting \(d\) and, in our experiments, generally chooses smaller values of \(d\) than the componentwise strategy, as illustrated in Figure~\ref{fig:gallery-experiments}.
If stronger componentwise stability is required, the additional cost of the componentwise strategy is justified by its tighter control of the local matrix-multiplication errors and the resulting componentwise backward error.



\section*{Acknowledgments}
Both authors are supported by the European Union (ERC, inEXASCALE, 101075632). Views and opinions expressed are those of the authors only and do not necessarily reflect those of the European Union or the European Research Council. Neither the European Union nor the granting authority can be held responsible for them. Both authors additionally acknowledge support from the Charles University Research Centre program No. UNCE/24/SCI/005.

\bibliographystyle{abbrvurl}
\bibliography{mybib}

\end{document}